\documentclass{CSML}

\def\dOi{11(4:10)2015}
\lmcsheading%
{\dOi}
{1--35}
{}
{}
{Mar.~11, 2015}
{Dec.~14, 2015}
{}

\ACMCCS{[{\bf Theory of computation}]: Models of computation---Computability---Turing machines}

\usepackage{hyperref}

\usepackage{amssymb}

\begin{document}

\title[Problems in number theory from busy beaver competition]{Problems in number theory\\ from busy beaver competition}
\author[P.~Michel]{Pascal MICHEL}
\address{\'Equipe de Logique Math\'ematique,
Institut de Math\'ematiques de Jussieu--Paris Rive Gauche, UMR 7586,
B\^atiment Sophie Germain, case 7012, 75205 Paris Cedex 13, France
and Universit\'e de Cergy-Pontoise, IUFM, F-95000 Cergy-Pontoise, France}
\email{michel@math.univ-paris-diderot.fr}
\thanks{Corresponding address: 59 rue du Cardinal Lemoine,
75005 Paris, France.}

\subjclass{Theory of computation, Models of computation, Computability, Turing machines}

\begin{abstract}
By introducing the busy beaver competition of Turing machines, in 1962, Rado defined
noncomputable functions on positive integers.
The study of these functions and variants leads to many mathematical challenges.
This article takes up the following one: How can a small Turing
machine manage to produce very big numbers?
It provides the following answer: mostly by simulating Collatz-like functions,
that are generalizations of the famous 3$x$+1 function.
These functions, like the 3$x$+1 function, lead to new unsolved problems in number theory.
\end{abstract}

\maketitle

\section{Introduction}\label{Sect1}
\subsection{A well defined noncomputable function}
It is easy to define a noncomputable function on nonnegative integers.
Indeed, given a programming language, you produce a systematic list of the programs for
functions, and, by diagonalization, you define a function whose output, on input $n$,
is different from the output of the $n$th program. This simple definition raises many problems:
Which programming language? How to list the programs? How to choose the output?

In 1962, Rado \cite{Ra62} gave a practical solution by defining the
\emph{busy beaver game}, also called now the \emph{busy beaver competition}.
Consider all Turing machines on one infinite tape, with $n$ states
(plus a special halting state), and two symbols (1, and the blank symbol 0),
and launch all of them on a blank tape. Define $S(n)$ as the maximum number
of computation steps made by such a machine before it stops, and define
$\Sigma(n)$ as the maximum number of symbols 1 left on the tape by a machine when it stops.
Then functions $S$ and $\Sigma$ are noncomputable, and, moreover,
grow faster than any computable function, that is, for any computable function $f$, there
exists an integer $N$ such that, for any $n \ge N$, $S(n) > \Sigma(n) > f(n)$.

More than fifty years later, no better choice has been found for a practical
noncomputable function. Only variants of Rado's definition have been proposed.
So, in 1988, Brady \cite{Br88} defined similar functions $S(n,m)$ and
$\Sigma(n,m)$ for $n \times m$ Turing machines, that is Turing machines with
$n$ states and $m$ symbols. He also introduced analogous functions for two-dimensional
Turing machines and ``turNing machines'', later resumed and expanded by Tim Hutton \cite{Hu}.
B\'atfai \cite{Ba1,Ba2} relaxed the rule about head moving, by allowing the head to stand still.

In this article, we will consider functions $S(n,m)$ and $\Sigma(n,m)$.
Recall that $S(n) = S(n,2)$ and $\Sigma(n) = \Sigma(n,2)$.

\subsection{Computing the values of noncomputable functions}
The busy beaver functions $S$ and $\Sigma$ are explicitely defined,
and it is possible to compute $S(n,m)$ and $\Sigma(n,m)$ for small $n$ and $m$.
In the first article on busy beavers, Rado \cite{Ra62} gave $\Sigma(2) = 2$ and
$\Sigma(3) \ge 6$. These results show that two problems are at stake:

\begin{itemize}
\item {\bf Problem 1}: To give lower bounds on $S(n,m)$ and $\Sigma(n,m)$
by finding Turing machines with high scores.
\item {\bf Problem 2}: To compute $S(n,m)$ and $\Sigma(n,m)$ by proving that
no Turing machines do better than the known best ones.
\end{itemize}

\noindent Problem 1 can be tackled either by hand search, as did, for example, Green \cite{Gr64} and Lynn \cite{Ly72},
or by computer search, using acceleration techniques of computation and, for example,
simulated annealing, as did T.\ and S.\ Ligocki \cite{Li}.

Solving Problem 2 requires more work to be done: clever enumeration of $n \times m$ Turing machines,
simulation of computation with acceleration techniques, proofs of non-halting for
the machines that do not halt.

Currently, the following results are known (see Michel \cite{MiA,MiS1} for a historical survey):
\begin{itemize}
\item $S(2) = 6$ and $\Sigma(2) = 4$ (Rado \cite{Ra62}),
\item $S(3) = 21$ and $\Sigma(3) = 6$ (Lin and Rado \cite{LR65}),
\item $S(4) = 107$ and $\Sigma(4) = 13$ (Brady \cite{Br66,Br83}, Machlin and Stout \cite{MS90}),
\item $S(5) \ge 47,176,870 $ and $\Sigma(5) \ge 4098$ (Marxen and Buntrock \cite{MB90}),
\item $S(6) > 7.4 \times 10^{36534}$ and $\Sigma(6) > 3.5 \times 10^{18267}$ (P.\ Kropitz in 2010),
\item $S(2,3) = 38$ and $\Sigma(2,3) = 9$ (Lafitte and Papazian \cite{LP07}),
\item $S(3,3) > 1.1 \times 10^{17}$ and $\Sigma(3,3) \ge 347,676,383$ (T.\ and S.\ Ligocki in 2007),
\item $S(4,3) > 1.0 \times 10^{14072}$ and $\Sigma(4,3) > 1.3 \times 10^{7036}$ (T.\ and S.\ Ligocki in 2008),
\item $S(2,4) \ge 3,932,964$ and $\Sigma(2,4) \ge 2050$ (T.\ and S.\ Ligocki in 2005),
\item $S(3,4) > 5.2 \times 10^{13036}$ and $\Sigma(3,4) > 3.7 \times 10^{6518}$ (T.\ and S.\ Ligocki in 2007),
\item $S(2,5) > 1.9 \times 10^{704}$ and $\Sigma(2,5) > 1.7 \times 10^{352}$ (T.\ and S.\ Ligocki in 2007),
\item $S(2,6) > 2.4 \times 10^{9866}$ and $\Sigma(2,6) > 1.9 \times 10^{4933}$ (T.\ and S.\ Ligocki in 2008).
\end{itemize}

In order to achieve these results, many computational and mathematical challenges had
to be taken up.

\begin{description}
\item[{\bf A}] Computational challenges.
\begin{enumerate}[label={\bf A\arabic*}.]
\item To generate all $n \times m$ Turing machines, or rather, to treat all
cases without having to generate all $n \times m$ Turing machines.
\item To simulate the computation of a machine by using acceleration
techniques (see Marxen and Buntrock \cite{MB90}, Marxen \cite{MaMa}).
\item To gave automatic proofs that non-halting machines do not halt
(see Brady \cite{Br83}, Marxen and Buntrock \cite{MB90}, Machlin and Stout \cite{MS90},
Hertel \cite{He09}, Lafitte and Papazian \cite{LP07}).
\end{enumerate}
\item[{\bf B}] Mathematical challenges.
\begin{enumerate}[label={\bf B\arabic*}.]
\item To prove by hand that a non-halting machine that resists the
computational proof does not halt.
\item To understand how the Turing machines that reach high scores manage to do it.
\end{enumerate}
\end{description}

\subsection{Facing open problems in number theory}
Let us come back to mathematical challenge \emph{B1}. For example,
the computational study of $5 \times 2$ Turing machines by Marxen and Buntrock
\cite{MB90}, Skelet \cite{Geo} and Hertel \cite{He09}
left holdouts that needed to be analyzed by hand. Marxen and Buntrock \cite{MB90} gave
an unsettled $5 \times 2$ Turing machine, named \#4, that turned out to never
halt, by an intricate analysis.

Actually, the halting problem for Turing machines launched on a blank tape
is m-complete, and this implies that this problem is as hard as the problem of the provability
of a mathematical statement in a logical theory such as ZFC (Zermelo Fraenkel set theory
with axiom of choice). So, when Turing machines with more and more states
and symbols are studied, potentially all theorems of mathematics will be met.
When more and more non-halting Turing machines are studied to be proved
non-halting, one has to expect to face hard open problems in mathematics,
that is problems that current mathematical knowledge can't settle.

\smallskip

Consider now mathematical challenge \emph{B2}, which is the very subject of
this article.

From 1983 to 1989, several $5 \times 2$ Turing machines with high scores were discovered
by Uwe Schult, by George Uhing, and by Heiner Marxen and J\"urgen Buntrock.
Michel \cite{MiT,Mi93} analyzed some of these machines and found that their behavior is
\emph{Collatz-like}, which implies that the halting problems on general inputs for these machines
are open problems in number theory (see Table \ref{52}).

From 2005 and 2007, many $3 \times 3$, $2 \times 4$ and $2 \times 5$ Turing machines with 
high scores were discovered, mainly by two teams: the French one of Gr\'egory Lafitte and Christophe Papazian,
and the father-and-son collaboration of Terry and Shawn Ligocki. Collatz-like behavior
of these champions seems to be the rule (see Tables \ref{33}, \ref{24} and \ref{25}).

However, the behaviors of $6 \times 2$ Turing machines display some variety. Many machines were
discovered, from 1990 to 2010, by Heiner Marxen and J\"urgen Buntrock, by Terry and Shawn Ligocki
and by Pavel Kropitz. The analyses of some of these machines, by Robert Munafo, Clive Tooth, Shawn
Ligocki and the author, show that the behaviors can be Collatz-like, exponential Collatz-like,
loosely Collatz-like, or definitely not Collatz-like. Almost all of them raise open problems
(see Table \ref{62}).

\emph{Note:} The Turing machines listed in Tables \ref{52}--\ref{25} are those for which
an analysis is known by the author. The machines without references for the study
of behavior were analyzed by the author \cite{MiA,MiS2}.
Many other machines are waiting for their analyses.

\subsection{Collatz functions, Collatz-like functions and other functions}
The 3$x$+1 \emph{function}, or \emph{Collatz function}, is the function $T$ on
positive integers defined by
$$T(n) = \left\{\begin{array}{cl}
n/2 & \mbox{if $n$ is even}\\
(3n+1)/2 & \mbox{if $n$ is odd}
\end{array}\right.$$
This function is famous because, when it is iterated on a positive integer, it seems
to lead to the loop $2,1,2,1,\ldots$. Is it always true? This is an open problem.
See Lagarias \cite{La85,LaB1,LaB2,La10} for more information.

It is natural to generalize the definition of the 3$x$+1 function by replacing $n$ even,
$n$ odd by $n \equiv 0,\dots,d-1\ (\mbox{mod}\ d)$, and by replacing $n/2$, $(3n+1)/2$ by $an+b$ for
rational numbers $a$, $b$. Unfortunately, no name for such functions is currently
taken for granted. Formal definitions were given
by Rawsthorne \cite{Ra85}, who proposed \emph{Collatz-type iteration functions},
by Buttsworth and Matthews \cite{BM90}, who proposed \emph{generalized Collatz mappings},
by Ka\v{s}\v{c}\'ak \cite{Ka92}, who proposed \emph{one-state linear operator algorithms (OLOA)},
and by Kohl \cite{Ko07}, who proposed \emph{residue-class-wise affine functions (RCWA)}.
Without giving a formal definition,
Lagarias \cite{La85} proposed \emph{periodically linear functions}, and Wagon \cite{Wa85}
proposed \emph{Collatz-like functions}.

We will choose the following definitions.
\begin{defi}\label{Def1.1}
A mapping $f:\mathbb{Z}\to\mathbb{Z}$ is a \emph{generalized Collatz mapping}
if there exists an integer $d \ge 2$ such that the following three equivalent conditions
are satisfied:
\begin{enumerate}[label=(\roman*)]
\item[(i)] (see \cite[p.14]{Wi98}) There exist rational numbers $q_0,\ldots,q_{d-1}$,
$r_0,\ldots,r_{d-1}$, such that, for all $i$, $0 \le i \le d-1$, we have
$q_id \in \mathbb{Z}$, $q_ii+r_i \in \mathbb{Z}$, and, for all $n \in \mathbb{Z}$,
$f(n) = q_in+r_i$ if $n \equiv i$ (mod $d$).
\item[(ii)] (see \cite{BM90}) There exist integers $m_0,\ldots,m_{d-1}$,
$p_0,\ldots,p_{d-1}$, such that, for all $i$, $0 \le i \le d-1$, we have
$p_i \equiv im_i$ (mod $d$) and, for all $n \in \mathbb{Z}$,
$f(n) = (m_in-p_i)/d$ if $n \equiv i$ (mod $d$).
\item[(iii)] There exist integers $a_0,\ldots,a_{d-1}$, $b_0,\ldots,b_{d-1}$,
such that we have, for all $i$, $0 \le i \le d-1$, for all $n \in \mathbb{Z}$,
$f(dn+i) = a_in + b_i$.
\end{enumerate}
\end{defi}

\noindent These definitions are easily seen to be equivalent: we have $a_i = m_i = q_id$
and $b_i = (im_i - p_i)/d = q_ii + r_i$.

The definitions above concern total functions, but, in this article, we always deal
with partial functions and functions with parameters, so we introduce the
following definitions.

\begin{defi}\label{Def1.2}
A partial function $f: \mathbb{Z} \to \mathbb{Z}$ is a
\emph{generalized Collatz function}, or a
\emph{pure Collatz-like function (without parameter)} if, in the previous definition,
$f(dn + i)$ can be undefined for one or many $i$, $0 \le i \le d-1$.
\end{defi}

\begin{defi}\label{Def1.3}
A partial function $f: \mathbb{Z} \times \mathbb{Z} \to \mathbb{Z} \times \mathbb{Z}$
is a \emph{pure Collatz-like function with parameter} if there exist an integer $d \ge 2$,
integers $a_0,\ldots,a_{d-1}$, $b_0,\ldots,b_{d-1}$, a set $S$ of integers and a function
$p : \{0,\ldots,d-1\} \times S \to S$ such that, for all $i$, $0 \le i \le d-1$, for all
$n \in \mathbb{Z}$, for all $s \in S$, $f(dn + i,s) = (a_in + b_i, p(i,s))$ or is undefined.
\end{defi}

\begin{defi}\label{Def1.4}
If, in the definitions above, $a_i = a$ for all $i$, $0 \le i \le d-1$,
we say that $f$ is pure Collatz-like \emph{of type} $d \to a$.
\end{defi}

We also need to define a new type of function, as follows.

\begin{defi}\label{Def1.5} A partial function $f : \mathbb{Z} \to \mathbb{Z}$ is an
\emph{exponential Collatz-like function} if there exist integers $d,p \ge 2$,
integers $a_0,\ldots,a_{d-1}$, $b_0,\ldots,b_{d-1}$,
$c_0,\ldots,c_{d-1}$, such that, for all $i$, $0 \le i \le d-1$,
all $n \in \mathbb{Z}$, $f(dn+i) = (a_ip^n + b_i)/c_i$ or is undefined.
In this definition, integers $p$, $a_i$, $b_i$, $c_i$ are chosen such that
$(a_ip^n + b_i)/c_i$ is an integer for all $n \in \mathbb{Z}$.
\end{defi}

Currently, no study of this type of function is known. Note that iterates
$f(n)$, $f^2(n),\ldots$ grow much faster for exponential Collatz-like functions
than for pure Collatz-like functions.

\subsection{From Collatz-like functions to high scores}

The Turing machines studied in this article have behaviors modeled on iterations
of functions, where halting configurations correspond to undefined values of functions.

How do Turing machines simulate Collatz-like functions?

First note that Baiocchi \cite{Ba98}, Margenstern \cite{Ma00} and Michel \cite{Mi14}
found Turing machines that simulate the 3$x$+1 function with a very small
number of states and symbols. In these articles, clever tricks were designed
to minimize the size of the Turing machines.

On the other hand, the Turing machines in the present article were frontally
attacked, and it is with handsight that they were found to simulate
Collatz-like functions.

Note also that, while the theorems giving the behavior of Turing machines
below are hard to be formally proven, they can be easily verified
for the small values of the parameters, by using programs simulating the Turing machines.

In Section \ref{Sect3}, we present a $3 \times 3$ Turing machine $M_1$ whose behavior is
pure Collatz-like, of type $8 \to 14$.
In Section \ref{Sect4}, we present a $2 \times 4$ Turing machine $M_2$ whose behavior is
pure Collatz-like with parameter, of type $3 \to 5$.
In Section \ref{Sect5}, we present a $2 \times 5$ Turing machine $M_3$ whose behavior is
pure Collatz-like with parameter, of type $2 \to 3$.
Thus, the halting problem for machines $M_1$, $M_2$ and $M_3$ depends on
an open problem about iterating Collatz-like functions.

In Section \ref{Sect6}, we present a $6 \times 2$ Turing machine $M_4$ whose behavior
is exponential Collatz-like.

In Section \ref{Sect7}, we present a $6 \times 2$ Turing machine $M_5$ whose behavior
depends on iterating a partial function $g_5(n,p)$. Without being Collatz-like,
this function seems to share some properties with Collatz-like functions.

In Section \ref{Sect8}, we present a $6 \times 2$ Turing machine $M_6$ whose behavior
looks like a loosely Collatz-like behavior with parameter, of type $2 \to 5$.
The novelty is that a potentially infinite set of rules seems to be necessary
to completely describe the behavior of the machine on inputs $00x$, $x \in \{0,1\}^*$.
A string $x \in \{0,1\}^*$ ending with symbol 1 can be taken as the binary
representation of a number $p$, read in the opposite direction, so $x = \mbox{R(bin($p$))}$,
where bin($p$) is the usual binary representation of $p$, and R($w$) is the reverse
of string $w$, that is R($w_1\ldots w_n$) = $w_n\ldots w_1$. In Table \ref{62}
we write ``R(bin($p$))'' to indicate the machines with a behavior involving
a potentially infinite set of rules. Of course, only a finite subset of these rules
are used when the machine is launched on a blank tape.

In Section \ref{Sect9}, we present a $6 \times 2$ Turing machine $M_7$ whose behavior
on the blank tape depends on configurations $C(n)$ all of them provably leading
to a halting configuration. We present such a machine to show how a Turing machine
can take a long time to stop without calling for Collatz-like functions.

\begin{table}[h]\centering
\begin{tabular}{|c|c|c|}
\hline
Machine            & Behavior                      & Study of behavior\\
\hline
January 1983       & Pure Collatz-like ($4\to 9$)  & Robinson,\\
Uwe Schult         & without parameter             & cited in \cite{De84}.\\
$\sigma = 501$     & 7 rules                       & Michel \cite{MiT} \\
$s$ = 134,467      & 7 transitions                 & \\
\hline
December 1984      & Pure Collatz-like ($3\to 8$)  & \\
George Uhing       & with parameter                & Michel \cite{MiT} \\
$\sigma = 1915$    & 5 rules                       & \\
$s$ = 2,133,492    & 9 transitions                 & \\
\hline
February 1986      & Pure Collatz-like ($8\to 15$) & \\
George Uhing       & without parameter             & \\
$\sigma = 1471$    & 5 rules                       &  \\
$s$ = 2,358,064    & 11 transitions                & \\
\hline
August 1989        & Pure Collatz-like ($3 \to 5$) & \\
Marxen, Buntrock   & without parameter             & Michel \cite{MiT} \\
$\sigma = 4098$    & 3 rules                       & \\
$s$ = 11,798,826   & 15 transitions                & \\
\hline
September 1989     & Pure Collatz-like ($3 \to 5$) & \\
Marxen, Buntrock   & without parameter             & Michel \cite{MiT} \\
$\sigma = 4097$    & 3 rules                       & \\
$s$ = 23,554,764   & 15 transitions                & \\
\hline
September 1989     & Pure Collatz-like ($3 \to 5$) & \\
Marxen, Buntrock   & without parameter             & Michel \cite{Mi93} \\
$\sigma = 4098$    & 3 rules                       & \\
$s$ = 47,176,870   & 15 transitions                & \\
\hline
\end{tabular}
\caption{\label{52}Study of behavior of $5\times 2$ machines.
For each machine, in the first column, one can find when it was discovered,
by whom, the number $\sigma$ of non-blank symbols left on the tape when the machine
halts, and the number $s$ of steps of the computation.
In the second column, the behavior of the machine is given, and we refer
to Def.\ 1.2--1.4 for the precise definition of pure Collatz-like function,
with and without parameter, of type $d \to a$. The number of rules gives roughly
the length of the definition of the function. The number of transitions
gives the number of times that the rules are used during the computation of the
machine on a blank tape. In the last column, machines without references were
analyzed by the author \cite{MiA,MiS2}.}
\end{table}

\begin{table}[h]\centering
\begin{tabular}{|c|c|c|}
\hline
Machine              & Behavior                      & Study of behavior\\
\hline
September 1997       & Pure Collatz-like ($4\to 10$) & \\
Marxen and Buntrock  & without parameter             & Munafo \cite{Mu1}\\
                     & 5 rules                       & \\
$s>8.6\times 10^{15}$ & 21 transitions                & \\
\hline
October 2000         &  R(bin(p)) ($2\to 3$)         & \\
Marxen and Buntrock  &                               & \\
(machine o)          & 9 rules                       & \\
$s>6.1\times 10^{119}$ & 337 transitions              & \\
\hline
October 2000          & All $C(n)$ stop             & \\
Marxen and Buntrock   &                              & Munafo \cite{Mu2}\\
(machine q)           & 4 rules                      & Section \ref{Sect9}\\
$s>6.1\times 10^{925}$ & 5 transitions                & \\
\hline
March 2001            & R(bin(p)) ($2\to 3$)         & \\
Marxen and Buntrock   &                              & Tooth \cite{To}\\
                      & 20 rules                     & \\
$s>3.0\times 10^{1730}$ & 4911 transitions            & \\
\hline
November 2007         & R(bin(p)) ($2\to 5$)         & \\
T. and S. Ligocki     &                              & Section \ref{Sect8}\\
                      & 12 rules                     & \\
$s>8.9\times 10^{1762}$ & 3346 transitions            & \\
\hline
December 2007         & R(bin(p)) ($4\to 6$)         & \\
T. and S. Ligocki     &                              & \\
                      & 18 rules                     & \\
$s>2.5\times 10^{2879}$ & 11026 transitions            & \\
\hline
May 2010               & Unclassifiable              & \\
Pavel Kropitz          &                             & S. Ligocki \cite{Li3}\\
                       & 6 rules                     & Section \ref{Sect7}\\
$s>3.8\times 10^{21132}$ & 22158 transitions           & \\
\hline
June 2010              & Exponential Collatz-like   & \\
Pavel Kropitz          & without parameter          & Section \ref{Sect6}\\
                       & 4 rules                    & \\
$s>7.4\times 10^{36534}$ & 5 transitions              & \\
\hline
\end{tabular}
\caption{\label{62}Study of behavior of $6\times 2$ machines.
See Def.\ 1.5 for the definition of exponential Collatz-like function.
We write R(bin(p)) when a potentially infinite set of rules would
be needed for a complete analysis of the machine.}
\end{table}

\begin{table}[h]\centering
\begin{tabular}{|c|c|c|}
\hline
Machine            & Behavior                      & Study of behavior\\
\hline
December 2004       & Pure Collatz-like ($2\to 5$) & \\
Brady               & with parameter               & \\
                    & 5 rules                      & \\
$s$ = 92,649,163    & 11 transitions               & \\
\hline
July 2005           & Pure Collatz-like ($2\to 5$) & \\
Souris              & with parameter               & \\
$\sigma = 36089$    & 5 rules                      & \\
$s$ = 310,341,163   & 12 transitions               & \\
\hline
July 2005           & Pure Collatz-like ($3\to 7$) & \\
Souris              & with parameter               & \\
                    & 7 rules                      & \\
$s$ = 544,884,219   & 12 transitions               & \\
\hline
August 2005         & Pure collatz-like ($4\to 7)$ & \\
Lafitte and Papazian & with parameter              & \\
                    & 8 rules                      & \\
$s>4.9\times 10^9$  & 21 transitions               & \\
\hline
September 2005      & Pure Collatz-like ($4\to 7$) & \\
Lafitte and Papazian & with parameter              & \\
                    & 7 rules                      & \\
$s>9.8\times 10^{11}$& 24 transitions               & \\
\hline
April 2006          & Pure Collatz-like ($2\to 5$) & \\
Lafitte and Papazian & with parameter              & \\
                    & 5 rules                      & \\
$s>4.1\times 10^{12}$& 16 transitions               & \\
\hline
August 2006         & Pure Collatz-like ($2\to 5$) & \\
T. and S. Ligocki   & with parameter               & S. Ligocki \cite{Li2}\\
                    & 4 rules                      & \\
$s>4.3\times 10^{15}$& 20 transitions               & \\
\hline
November 2007       & Pure Collatz-like ($8\to 14$)& \\
T. and S. Ligocki   & without parameter            & Section \ref{Sect3}\\
                    & 9 rules                      & \\
$s>1.1\times 10^{17}$& 34 transitions               & \\
\hline
\end{tabular}
\caption{\label{33}Study of behavior of $3\times 3$ machines}
\end{table}

\begin{table}[h]\centering
\begin{tabular}{|c|c|c|}
\hline
Machine           & Behavior                     & Study of behavior\\
\hline
1988              & Pure Collatz-like ($3\to 5$) & \\
Brady             & with parameter               & \\
                  & 6 rules                      & \\
$s$ = 7195        & 7 transitions                & \\
\hline
February 2005     & Pure Collatz-like ($3\to 5$) & \\
T. and S. Ligocki & with parameter               & Section \ref{Sect4}\\
                  & 7 rules                      & \\
$s$ = 3,932,964   & 14 transitions               & \\
\hline
\end{tabular}
\caption{\label{24}Study of behavior of $2\times 4$ machines}
\end{table}

\begin{table}[h]\centering
\begin{tabular}{|c|c|c|}
\hline
Machine               & Behavior                     & Study of behavior\\
\hline
October 2005          & Pure Collatz-like ($2\to 5$) & \\
Lafitte, Papazian     & with parameter               & \\
                      & 7 rules                      & \\
$s>9.1\times 10^{11}$  & 15 transitions               & \\
\hline
December 2005         & Pure Collatz-like ($2\to 5$) & \\
Lafitte, Papazian     & with parameter               & \\
                      & 5 rules                      & \\
$s>9.2\times 10^{11}$  & 14 transitions               & \\
\hline
May 2006              & Pure Collatz-like ($3\to 4$) & \\
Lafitte, Papazian     & with parameter               & \\
                      & 7 rules                      & \\
$s>3.7\times 10^{12}$  & 45 transitions               & \\
\hline
June 2006             & Pure Collatz-like ($2\to 3$) & \\
Lafitte, Papazian     & with parameter               & \\
                      & 9 rules                      & \\
$s>1.4\times 10^{13}$  & 36 transitions               & \\
\hline
August 2006           & Pure Collatz-like ($2\to 5$) & \\
T. and S. Ligocki     & with parameter               & S. Ligocki \cite{Li1}\\
                      & 9 rules                      & \\
$s>7.0\times 10^{21}$  & 30 transitions               & \\
\hline
November 2007         & Pure Collatz-like ($2\to 3$)   & \\
T. and S. Ligocki     & with parameter               & Section \ref{Sect5}\\
                      & 17 rules                     & \\
$s>1.9\times 10^{704}$ & 2002 transitions             & \\
\hline
\end{tabular}
\caption{\label{25}Study of behavior of $2\times 5$ machines}
\end{table}

\section{Preliminaries}
A \emph{Turing machine} involved in the busy beaver competition is defined as follows.
It has a tape made of cells, infinite in both directions.
Each cell contains a \emph{symbol}, and one head can move on the tape and read and write a symbol on a cell.
The Turing machine can be in a finite number of \emph{states}.
A computation of the Turing machine is a sequence of steps.
In a step of computation, 
according to the current state and the symbol read by the head on the current cell,
the head writes a symbol on the cell, moves to the next cell on the right side
or on the left side, and the machine enters a new state. 

\noindent Formally, a Turing machine $M = (Q,\Gamma,\delta)$ has a finite set of states
$Q = \{q_0,q_1,\ldots,q_{n-1}\}$, a finite set of symbols $\Gamma = \{0,1,\ldots,m-1\}$, and
a \emph{transition function} (or \emph{next move function}) $\delta$,
which is a mapping
$$\delta : Q \times \Gamma \to (\Gamma \times \{L,R\} \times Q)\cup\{(1,R,H)\}.$$
If $\delta(q,a) = (b,D,q')$, then the Turing machine,
when it is in state $q$ reading symbol $a$ on the current cell, 
writes symbol $b$ instead of $a$ on this cell, moves one cell left if $D = L$,
one cell right if $D = R$, and enters state $q'$.
The transition function is usually given by a \emph{transition table}. 

 For greater convenience, the states of the Turing machines will be denoted by capital
letters: $A$, $B,\ldots$.
There is a special state $A$, called the \emph{initial state}, and a special symbol
$0$, called the \emph{blank symbol}. In the busy beaver competition, at the beginning of
a computation, the Turing machine is in state $A$, and the tape is blank,
that is all the cells of the tape contain the blank symbol.
There is another state $H$, the \emph{halting state}, not in the set $Q$ of states.
When a Turing machine enters this state, the computation stops.
We impose that, at the last step, the machine writes 1, moves right,
and enters state $H$. The machines are \emph{normalized} in the following way:
When they are launched on a blank tape, they enter new states in the order $B$,
$C,\ldots$, and they write new symbols in the order $1, 2,\ldots$. 

A \emph{word} is a finite string of symbols. The set of words with symbols in
the set $\Gamma$ is denoted by $\Gamma^*$. The number of symbols in a word
$x \in \Gamma^*$ is called the \emph{length} of $x$ and is denoted by
$|x|$. The \emph{empty word} is the word of length zero, denoted by $\lambda$.
If $x \in \Gamma^*$, and $n \ge 0$, $x^n$ is the word $xx\ldots x$, where $x$ is repeated
$n$ times, that is, formally: $x^0 = \lambda$, $x^1 = x$ and $x^{n+1} = x^nx$.
An infinite to the left string of 0 is denoted by $^\omega 0$, and an infinite to the
right string of 0 is denoted by $0^\omega$.

A \emph{configuration} is a way to encode the symbols on the tape, the state,
and the cell currently read by the head. The Turing machine is in configuration
$^\omega 0x(Sa)y0^\omega$, with $S \in Q\cup\{H\}$, $a \in \Gamma$,
$x,y \in \Gamma^*$, if the word $xay$ is written on the tape, the state is $S$, and
the head is reading symbol $a$.
Since, at the beginning of the computation, the state is $A$ and the tape is blank,
the initial configuration is $^\omega0(A0)0^\omega$.
If the state is $H$, the configuration is halting.
We also consider configurations $x(Sa)y$ with finite length.
If the computation from configuration $C_1$ to configuration $C_2$ takes $t$ steps, we write
$C_1 \vdash(t)\ C_2$, and $t$ is said to be the \emph{time} taken by the machine
to go from $C_1$ to $C_2$.
If $C_2$ is a halting configuration, we also write $C_1 \vdash(t)$ END.
We write $C_1 \vdash(\ )\ C_2$ if the time is not specified.
If $C_1$ and $C_2$ are configurations with finite length, then they refer to the
same part of the tape. For example, $(A0)0\vdash(1)\ 1(B0)$ if $\delta(A,0) = (1,R,B)$.

A Turing machine $M$ computes a partial function $f_M: \Gamma^* \to \Gamma^*$ as follows.
let $x = x_1\ldots x_n \in \Gamma^*$. Then $x$ becomes an input for $M$ by considering
the computation of $M$ on initial configuration $^\omega0(Ax_1)x_2\ldots x_n0^\omega$.
If $M$ never stops on this configuration, then $f_M(x)$ is undefined. If $M$ stops,
in configuration $^\omega0y(Ha)z0^\omega$, with $a \in \Gamma$, $y,z \in \Gamma^*$,
then the output $f_M(x)$ is defined from this configuration by a suitable convention.
The \emph{halting set} is $\{x \in \Gamma^*: f_M(x) \mbox{\ is defined}\}$.
The \emph{halting problem} for machine $M$ is the problem consisting in determining
the halting set. Note that the Turing machines with two symbols 0 and 1 are powerful
enough to compute any computable function, and their halting sets can be any
computably enumerable (also called recursively enumerable) set.

A Turing machine with $n$ states and $m$ symbols is called a $n \times m$ machine.
The set of $n \times m$ machines is denoted by TM($n$,$m$).
With our definition of the transition function, there are $(2nm+1)^{nm}$
machines in the set TM($n$,$m$). In the \emph{busy beaver competition}, for fixed
numbers of states $n$ and symbols $m$, all the $(2nm+1)^{nm}$ Turing machines in
TM($n$,$m$) are launched on the blank tape. Some of them never stop. Those which
stop are called \emph{busy beaver}. Each busy beaver takes some time to stop,
and leaves some non-blank symbols on the tape, so busy beavers are involved
in two competitions: to take the longest time before stopping, and to leave
the greatest number of non-blank symbols on the tape when stopping.
The time taken by Turing machine $M$ to stop is denoted by $s(M)$,
and the number of non-blank symbols left by $M$ when it stops is denoted
by $\sigma(M)$. The \emph{busy beaver functions} are defined by
$$S(n,m) = \max\{s(M) : M\ \mbox{is a busy beaver with $n$ states and $m$ symbols}\}$$
$$\Sigma(n,m) = \max\{\sigma(M) : M\ \mbox{is a busy beaver with $n$ states and $m$ symbols}\}$$
Rado \cite{Ra62} initially defined functions $S(n) = S(n,2)$ and $\Sigma(n) = \Sigma(n,2)$
for Turing machines with $n$ states and two symbols.

\section{Pure Collatz-like behavior}\label{Sect3}

Let $M_1$ be the $3 \times 3$ Turing machine defined by Table~\ref{m1}

\begin{table}[t]\centering
\begin{tabular}{|c|c|c|c|}
\hline
$M_1$ &  0  &  1  &  2\\
\hline
$A$   & 1RB & 2LA & 1LC\\
\hline
$B$   & 0LA & 2RB & 1LB\\
\hline
$C$   & 1RH & 1RA & 1RC\\
\hline
\end{tabular}
\caption{\label{m1}Machine $M_1$ discovered in November 2007 by T.\ and S.\ Ligocki.
Such tables are read as in the following example. When machine $M_1$ is in state $A$
and reads symbol 0, then it writes symbol 1 instead of symbol 0, moves one cell to the right,
and enters state $B$.}
\end{table}

We have $s(M_1)$ = 119,112,334,170,342,540 and $\sigma(M_1)$ = 374,676,383.

This machine is the current champion for the busy beaver competition for $3\times 3$
machines. It was discovered in November 2007 by Terry and Shawn Ligocki, who wrote
(email on November, 9th) that they enumerated all the $3 \times 3$ machines
and applied the techniques of acceleration and proof systems originally
developed by Marxen and Buntrock.

The following theorem gives the rules that enable Turing machine $M_1$ to reach
a halting configuration from a blank tape.

\begin{thm}\label{Thm3.1}
We have the following transitions between configurations of Turing machine $M_1$.
Let $C(n)$ = $^\omega 0(A0)2^n0^\omega$. Then

$\begin{array}{clll}
(a) & ^\omega 0(A0)0^\omega & \vdash(3) &  C(1),
\end{array}$

and, for all $k \ge 0$,

$\begin{array}{clll}
(b) & C(8k + 1) & \vdash(112k^2 + 116k + 13)  & C(14k + 3),\\
(c) & C(8k + 2) & \vdash(112k^2 + 144k + 38)  & C(14k + 7),\\
(d) & C(8k + 3) & \vdash(112k^2 + 172k + 54)  & C(14k + 8),\\
(e) & C(8k + 4) & \vdash(112k^2 + 200k + 74)  & C(14k + 9),\\
(f) & C(8k + 5) & \vdash(112k^2 + 228k + 97)  & ^\omega 01(H1)2^{14k + 9}0^\omega,\\
(g) & C(8k + 6) & \vdash(112k^2 + 256k + 139) & C(14k + 14),\\
(h) & C(8k + 7) & \vdash(112k^2 + 284k + 169) & C(14k + 15),\\
(i) & C(8k + 8) & \vdash(112k^2 + 312k + 203) & C(14k + 16).
\end{array}$
\end{thm}

\proof
A direct inspection of the transition table gives

$\begin{array}{clll}
(1)  &  0(A0)0   & \vdash(3)  & (A0)20,\\
(2)  & 0^3(A0)2^5 & \vdash(53) & (B1)1^8,\\
(3)  & 0(A1)     & \vdash(1)  & (A0)2,\\
(4)  & 1(A1)     & \vdash(1)  & (A1)2,\\
(5)  & 02(A1)    & \vdash(3)  & 1(H1)2,\\
(6)  & 12(A1)    & \vdash(4)  & (A1)22,\\
(7)  & 22(A1)    & \vdash(8)  & (A1)22,\\
(8)  & 2(B1)0^2  & \vdash(7)  & 11(A1)0,\\
(9)  & (B1)1     & \vdash(1)  & 2(B1),\\
(10) & (B1)2     & \vdash(1)  & 2(B2),\\
(11) & 0^3(B2)   & \vdash(14) & 1^3(B1),\\
(12) & 1(B2)     & \vdash(1)  & (B1)1,\\
(13) & 2(B2)     & \vdash(1)  & (B2)1.
\end{array}$

\smallskip

\noindent From this point, $k$ will be an integer, $k \ge 0$.

\smallskip

\noindent Iterating, respectively, (4), (7), (9) and (13) gives

$\begin{array}{clll}
(14) & 1^k(A1)   & \vdash(k)  & (A1)2^k,\\
(15) & 2^{2k}(A1) & \vdash(8k) & (A1)2^{2k},\\
(16) & (B1)1^k   & \vdash(k)  & 2^k(B1),\\
(17) & 2^k(B2)   & \vdash(k)  & (B2)1^k.
\end{array}$

\smallskip

\noindent Using consecutively (16), (10), (17) and (12), we get
 
$\begin{array}{clll}
(18) & 1(B1)1^k2 & \vdash(2k+3) & (B1)1^{k+2}.
\end{array}$

\smallskip

\noindent Using (16), (10), (17) and (11), we get

$\begin{array}{clll}
(19) & 0^3(B1)1^k 2& \vdash(2k+16) & 1^3(B1)1^{k+1}.
\end{array}$

\smallskip

\noindent Using (19) and three times (18), we get

$\begin{array}{clll}
(20) & 0^3(B1)1^k2^4 & \vdash(8k+43) & (B1)1^{k+7}.
\end{array}$

\smallskip

\noindent For any $n \ge 0$, by induction on $k$, using (20), we get

$\begin{array}{clll}
(21) & 0^{3k}(B1)1^n2^{4k} & \vdash(28k^2+(8n+15)k) & (B1)1^{7k+n}.
\end{array}$

\smallskip

\noindent By taking $n=8$ in (21), we get

$\begin{array}{clll}
(22) & 0^{3k}(B1)1^82^{4k} & \vdash(28k^2+79k) & (B1)1^{7k+8}.
\end{array}$

\smallskip

\noindent Using (8), (14) and (15), we get

$\begin{array}{clll}
(23) & 2^{2k+1}(B1)0^2 & \vdash(8k+9) & (A1)2^{2k+2}0.
\end{array}$

\smallskip

We are now ready to prove the results of the theorem.

\noindent Using (2), (22) and (16), we get

$\begin{array}{clll}
(24) & 0^{3k+3}(A0)2^{4k+5} & \vdash(28k^2 + 86k + 61) & 2^{7k+8}(B1).
\end{array}$

\smallskip

\noindent Using (24), (23) and (5) we get

$\begin{array}{clll}
(25) & 0^{6k+4}(A0)2^{8k+5}0^2 & \vdash(112k^2 + 228k + 97) & 1(H1)2^{14k+9}0.
\end{array}$

\smallskip

\noindent Using (24), (23) and (3) we get

$\begin{array}{clll}
 & 0^{6k+7}(A0)2^{8k+9}0^2 & \vdash(112k^2 + 340k + 241) & (A0)2^{14k+17}0,
\end{array}$

\noindent and this result is still true for $k = -1$, so we have

$\begin{array}{clll}
(26) & 0^{6k+1}(A0)2^{8k+1}0^2 & \vdash(112k^2 + 116k + 13) & (A0)2^{14k+3}0.
\end{array}$

\smallskip

\noindent Using (2), (22), (19) and (16) we get

$\begin{array}{clll}
(27) & 0^{3k+6}(A0)2^{4k+6} & \vdash(28k^2 + 100k + 94) & 1^32^{7k+9}(B1).
\end{array}$

\smallskip

\noindent Using (27), (23), (14) and (3) we get

$\begin{array}{clll}
(28) & 0^{6k+7}(A0)2^{8k+6}0^2 & \vdash(112k^2 + 256k + 139) & (A0)2^{14k+14}0.
\end{array}$

\smallskip

\noindent Using (27), (23), (6), (14) and (3) we get

$\begin{array}{clll}
 & 0^{6k+10}(A0)2^{8k+10}0^2 & \vdash(112k^2 + 368k + 294) & (A0)2^{14k+21}0,
\end{array}$

\noindent and this result is still true for $k = -1$, so we have

$\begin{array}{clll}
(29) & 0^{6k+4}(A0)2^{8k+2}0^2 & \vdash(112k^2 + 144k + 38) & (A0)2^{14k+7}0.
\end{array}$

\smallskip

\noindent Using (2), (22), (19), (18), (16), (8) and (14) we get

$\begin{array}{clll}
(30) & 0^{3k+6}(A0)2^{4k+7}0^2 & \vdash(28k^2 + 114k + 126) & 1^22^{7k+10}(A1)220.
\end{array}$

\smallskip

\noindent Using (30), (15), (14) and (3) we get

$\begin{array}{clll}
(31) & 0^{6k+7}(A0)2^{8k+7}0^2 & \vdash(112k^2 + 284k + 169) & (A0)2^{14k+15}0.
\end{array}$

\smallskip

\noindent Using (30), (15), (6), (4) and (3) we get

$\begin{array}{clll}
 & 0^{6k+10}(A0)2^{8k+11}0^2 & \vdash(112k^2 + 396k + 338) & (A0)2^{14k+22}0,
\end{array}$

\noindent and this result is still true for $k = -1$, so we have

$\begin{array}{clll}
(32) & 0^{6k+4}(A0)2^{8k+3}0^2 & \vdash(112k^2 + 172k + 54) & (A0)2^{14k+8}0.
\end{array}$

\smallskip

\noindent Using (2), (22), (19), (18), (18), (16), (8) and (14) we get

$\begin{array}{clll}
(33) & 0^{3k+6}(A0)2^{4k+8}0^2 & \vdash(28k^2 + 128k + 153) & 12^{7k+12}(A1)220.
\end{array}$

\smallskip

\noindent Using (33), (15), (4) and (3) we get

$\begin{array}{clll}
(34) & 0^{6k+7}(A0)2^{8k+8}0^2 & \vdash(112k^2 + 312k + 203) & (A0)2^{14k+16}0.
\end{array}$

\smallskip

\noindent Using (33), (15), (6) and (3) we get

$\begin{array}{clll}
 & 0^{6k+10}(A0)2^{8k+12}0^2 & \vdash(112k^2 + 424k + 386) & (A0)2^{14k+23}0,
\end{array}$

\noindent and this result is still true for $k = -1$, so we have

$\begin{array}{clll}
(35) & 0^{6k+4}(A0)2^{8k+4}0^2 & \vdash(112k^2 + 200k + 74) & (A0)2^{14k+9}0.
\end{array}$

\smallskip

\noindent The results (1), (26), (29), (32), (35), (25), (28), (31) and (34)
give, respectively, the results (a)--(i) of the theorem. \qed 

\medskip

\noindent Using the rules of this theorem, we have, in 34 transitions,
$$^\omega0(A0)0^\omega \vdash(3)\quad C(1) \vdash(13)\quad C(3) \vdash(\ ) \cdots
\vdash(\ )\quad ^\omega01(H1)2^{374,676,381}0^\omega.$$

\noindent Let us try to give an informal explanation of how Turing machine $M_1$ works.
Note that this explanation will also apply to the machines $M_2$ and $M_3$
below. The two main ideas are the following ones.
First, a string of symbols 2 is nibbled four by four, and each time a string of
four symbols 2 is nibbled, a string of seven symbols 2 is potentially created
(see transitions (20) and (21) in the proof of Theorem \ref{Thm3.1}).
This explains the type $8 \to 14$ of the simulated Collatz-like function.
Second, there are edge effects when the ends of a string of symbols 2 are reached,
which explain the multiplicity of cases.

Let $g_1$ be the pure Collatz-like function defined by: for $k \ge 0$,

$$\begin{array}{rcl}
g_1(8k + 1) & = & 14k + 3,\\
g_1(8k + 2) & = & 14k + 7,\\
g_1(8k + 3) & = & 14k + 8,\\
g_1(8k + 4) & = & 14k + 9,\\
g_1(8k + 5) &   & \mbox{undefined},\\
g_1(8k + 6) & = & 14k + 14,\\
g_1(8k + 7) & = & 14k + 15,\\
g_1(8k + 8) & = & 14k + 16.
\end{array}$$
  
\smallskip

\noindent Then $g_1^{33}(1)$ is undefined.

The theorem gives immediately the following proposition.

\begin{prop}\label{Prop3.2}
The behavior of Turing machine $M_1$, on inputs $02^n$, $n \ge 1$,
depends on the behavior of iterated $g_1^k(n)$, $k \ge 1$.\qed
\end{prop}

\noindent Since the behavior of iterated $g_1^k(n)$ is an open problem in mathematics,
this is also the case for the halting problem for Turing machine $M_1$.

Let $h_1(n) = \min\{k:g_1^k(n) \mbox{\ is undefined}\}$. We have seen that
$h_1(1) = 33$. We also have $h_1(144) = 41$, $h_1(270) = 51$.

\section{Collatz-like with parameter: first example}\label{Sect4}

Let $M_2$ be the $2 \times 4$ Turing machine defined by Table~\ref{m2}.

\begin{table}[t]\centering
\begin{tabular}{|c|c|c|c|c|}
\hline
$M_2$ &  0  &  1  &  2  & 3\\
\hline
$A$   & 1RB & 2LA & 1RA & 1RA\\
\hline
$B$   & 1LB & 1LA & 3RB & 1RH\\
\hline
\end{tabular}
\caption{\label{m2}Machine $M_2$ discovered in February 2005 by T.\ and S.\ Ligocki}
\end{table}

We have $s(M_2)$ = 3,932,964 and $\sigma(M_2)$ = 2050.

This machine is the current champion for the busy beaver competition for $2\times 4$
machines. It was discovered in February 2005 by Terry and Shawn Ligocki, who wrote
(email on February, 13th) that they found this machine using simulated annealing.

The following theorem gives the rules that enable Turing machine $M_2$ to reach
a halting configuration from a blank tape.

\begin{thm}\label{Thm4.1}
We have the following transitions between configurations of Turing machine $M_2$.
Let

$\begin{array}{lcl}
C(n,1) & = & ^\omega 0(A0)2^n10^\omega,\\
C(n,2) & = & ^\omega 0(A0)2^n110^\omega.
\end{array}$

Then

$\begin{array}{clll}
(a) & ^\omega 0(A0)0^\omega & \vdash(6) &  C(1,2),
\end{array}$

and, for all $k \ge 0$,

$\begin{array}{clll}
(b) &  C(3k,1)    & \vdash(15k^2 + 9k + 3)   & C(5k + 1,1),\\
(c) & C(3k + 1,1) & \vdash(15k^2 + 24k + 13) & ^\omega 013^{5k+2}1(H1)0^\omega,\\
(d) & C(3k + 2,1) & \vdash(15k^2 + 29k + 17) & C(5k + 4,2),\\
(e) &  C(3k,2)    & \vdash(15k^2 + 11k + 3)  & C(5k + 1,2),\\
(f) & C(3k + 1,2) & \vdash(15k^2 + 21k + 7)  & C(5k + 3,1),\\
(g) & C(3k + 2,2) & \vdash(15k^2 + 36k + 23) & ^\omega 013^{5k+4}1(H1)0^\omega.
\end{array}$
\end{thm}

\proof
A direct inspection of the transition table gives

$\begin{array}{clll}
(1)  & 0^2(A0)0 & \vdash(6)  & (A0)211,\\
(2)  & 1(A0)0   & \vdash(3) & (A1)11,\\
(3)  & (A0)2    & \vdash(1) & 1(B2),\\
(4)  & 0(A1)    & \vdash(1) & (A0)2,\\
(5)  & 1(A1)    & \vdash(1) & (A1)2,\\
(6)  & 3(A1)    & \vdash(2) & 1(A2),\\
(7)  & (A2)0    & \vdash(1) & 1(A0),\\
(8)  & (A2)1    & \vdash(1) & 1(A1),\\
(9)  & (A2)2    & \vdash(1) & 1(A2),\\
(10) & (B2)0    & \vdash(3) & 1(H1),\\
(11) & (B2)1    & \vdash(4) & (A1)2,\\
(12) & (B2)2    & \vdash(1) & 3(B2).
\end{array}$

\smallskip

\noindent Iterating, respectively, (5), (9) and (12) gives

$\begin{array}{clll}
(13) & 1^k(A1)   & \vdash(k)  & (A1)2^k,\\
(14) & (A2)2^k & \vdash(k) & 1^k(A2),\\
(15) & (B2)2^k   & \vdash(k)  & 3^k(B2).\\
\end{array}$

\smallskip

\noindent Using (2), (13) and (4) we get

$\begin{array}{clll}
 & 01^{k+1}(A0)0 & \vdash(k + 4) & (A0)2^{k+1}11,
\end{array}$

\noindent and this result is still true for $k = -1$, so we have

$\begin{array}{clll}
(16) & 01^k(A0)0 & \vdash(k + 3) & (A0)2^k11.
\end{array}$

\smallskip

\noindent Using (2), (13), (6), (14), (8), (13) and (4) we get

$\begin{array}{clll}
(17) & 0131^{k+1}(A0)0 & \vdash(3k + 10) & (A0)2^{k+4}1.
\end{array}$

\smallskip

\noindent Using (2), (13), (6), (14), (8), (13), (6), (14), (8), (13) and (4) we get

$\begin{array}{clll}
(18) & 01331^{k+1}(A0)0 & \vdash(5k+19) & (A0)2^{k+6}.
\end{array}$

\smallskip

\noindent Using (3), (15) and (10) we get

$\begin{array}{clll}
(19) & (A0)2^{k+1}0 & \vdash(k+4) & 13^k1(H1).
\end{array}$

\smallskip

\noindent Using (3), (15), (11), (6), (9) and (7) we get

$\begin{array}{clll}
(20) & (A0)2^{k+2}10 & \vdash(k+10) & 13^k1^3(A0).
\end{array}$

\smallskip

\noindent Using (3), (15), (11), (6), (9), (8), (13), (6), (14) and (7) we get

$\begin{array}{clll}
(21) & (A0)2^{k+3}110 & \vdash(k+20) & 13^k1^5(A0).
\end{array}$

\smallskip

\noindent Using (2), (13), (6), (14), (8), (13), (6), (14), (8), (13), (6), (14) and (7) we get

$\begin{array}{clll}
(22) & 3^31^{k+1}(A0)0^2 & \vdash(6k+24) & 1^{k+6}(A0).
\end{array}$

\smallskip

\noindent By induction on $k$, from (22), we get

$\begin{array}{clll}
(23) & 3^{3k}1^{n+1}(A0)0^{2k} & \vdash(15k^2 + 6nk + 9k) & 1^{5k+n+1}(A0),
\end{array}$

\noindent so we have, for $n = 2$ and $n = 4$

$\begin{array}{clll}
(24) & 3^{3k}1^3(A0)0^{2k} & \vdash(15k^2 + 21k) & 1^{5k+3}(A0),\\
(25) & 3^{3k}1^5(A0)0^{2k} & \vdash(15k^2 + 33k) & 1^{5k+5}(A0).
\end{array}$

\smallskip

We are now ready to prove the theorem.

\noindent Using (20), (24) and (17) we get

$\begin{array}{clll}
 & 0(A0)2^{3k+3}10^{2k+2} & \vdash(15k^2 + 39k + 27) & (A0)2^{5k+6}1,
\end{array}$

\noindent and the result is still true for $k = -1$, so we have

$\begin{array}{clll}
(26) & 0(A0)2^{3k}10^{2k} & \vdash(15k^2 + 9k + 3) & (A0)2^{5k+1}1.
\end{array}$

\smallskip

\noindent Using (20), (24), (18) and (19) we get

$\begin{array}{clll}
 & 0(A0)2^{3k+4}10^{2k+3} & \vdash(15k^2 + 54k + 52) & 13^{5k+7}1(H1),
\end{array}$

\noindent and the result is still true for $k = -1$, so we have

$\begin{array}{clll}
(27) & 0(A0)2^{3k+1}10^{2k+1} & \vdash(15k^2 + 24k + 13) & 13^{5k+2}1(H1).
\end{array}$

\smallskip

\noindent Using (20), (24) and (16) we get

$\begin{array}{clll}
(28) & 0(A0)2^{3k+2}10^{2k+2} & \vdash(15k^2 + 29k + 17) & (A0)2^{5k+4}11.
\end{array}$

\smallskip

\noindent Using (21), (25) and (16) we get

$\begin{array}{clll}
 & 0(A0)2^{3k+3}110^{2k+2} & \vdash(15k^2 + 41k + 29) & (A0)2^{5k+6}11,
\end{array}$

\noindent and the result is still true for $k = -1$, so we have

$\begin{array}{clll}
(29) & 0(A0)2^{3k}110^{2k} & \vdash(15k^2 + 11k + 3) & (A0)2^{5k+1}11.
\end{array}$

\smallskip

\noindent Using (21), (25) and (17) we get

$\begin{array}{clll}
 & 0(A0)2^{3k+4}110^{2k+2} & \vdash(15k^2 + 51k + 43) & (A0)2^{5k+8}1,
\end{array}$

\noindent and the result is still true for $k = -1$, so we have

$\begin{array}{clll}
(30) & 0(A0)2^{3k+1}110^{2k} & \vdash(15k^2 + 21k + 7) & (A0)2^{5k+3}1.
\end{array}$

\smallskip

\noindent Using (21), (25), (18) and (19) we get

$\begin{array}{clll}
 & 0(A0)2^{3k+5}110^{2k+3} & \vdash(15k^2 + 66k + 74) & 13^{5k+9}1(H1),
\end{array}$

\noindent and the result is still true for $k = -1$, so we have

$\begin{array}{clll}
(31) & 0(A0)2^{3k+2}110^{2k+1} & \vdash(15k^2 + 36k + 23) & 13^{5k+4}1(H1).
\end{array}$

\smallskip

\noindent The theorem comes from results (1) and (26)--(31). \qed

\smallskip

\noindent Using the rules of this theorem, we have, in 14 transitions,
$$^\omega0(A0)0^\omega \vdash(6)\quad C(1,2) \vdash(7)\quad C(3,1) \vdash(\ ) \cdots
\vdash(\ )\quad ^\omega013^{2047}1(H1)0^\omega.$$
Informally, Turing machine $M_2$ works according to the same ideas as Turing
machine $M_1$ does, but, in addition, edge effects depend also on the parameter.

Let $g_2$ be the pure Collatz-like function with parameter defined by: for $k \ge 0$,

$$\begin{array}{rcl}
g_2(3k, 1) & = & (5k + 1, 1),\\
g_2(3k + 1, 1) &   & \mbox{undefined},\\
g_2(3k + 2, 1) & = & (5k + 4, 2),\\
g_2(3k, 2) & = & (5k + 1, 2),\\
g_2(3k + 1, 2) & = & (5k + 3, 1),\\
g_2(3k + 2, 2) &   & \mbox{undefined}.
\end{array}$$
Then $g_2^{13}(1,2)$ is undefined.

The theorem gives immediately the following proposition.

\begin{prop}\label{Prop4.2}
The behavior of Turing machine $M_2$, on inputs $02^n1^i$, $n \ge 1$, $i \in \{1,2\}$
depends on the behavior of iterated $g_2^k(n,i)$, $k \ge 1$.\qed
\end{prop}

\noindent Since the behavior of iterated $g_2^k(n,i)$ is an open problem in mathematics,
this is also the case for the halting problem for Turing machine $M_2$.

Let $h_2(n,i) = \min\{k:g_2^k(n,i) \mbox{\ is undefined}\}$. We have seen that
$h_2(1,2) = 13$. We also have $h_2(137,1) = 16$, $h_2(210,2) = 20$.

\section{Collatz-like with parameter: second example}\label{Sect5}

Let $M_3$ be the $2 \times 5$ Turing machine defined by Table~\ref{m3}.

\begin{table}[t]\centering
\begin{tabular}{|c|c|c|c|c|c|}
\hline
$M_3$ &  0  &  1  &  2  &  3  &  4\\
\hline
$A$   & 1RB & 2LA & 1RA & 2LB & 2LA\\
\hline
$B$   & 0LA & 2RB & 3RB & 4RA & 1RH\\
\hline
\end{tabular}
\caption{\label{m3}Machine $M_3$ discovered in November 2007 by T.\ and S.\ Ligocki}
\end{table}

We have $s(M_3) > 1.9 \times 10^{704}$ and $\sigma(M_3) > 1.7 \times 10^{352}$.

This machine is the current champion for the busy beaver competition for $2 \times 5$
machines. It was discovered in November 2007 by Terry and Shawn Ligocki, who wrote
(email on November, 9th) that, as they did for $3 \times 3$ machine $M_1$,
they enumerated all the $2 \times 5$ machines and applied the techniques
of acceleration and proof systems originally developed by Marxen and Buntrock.

The following theorem gives the rules that enable Turing machine $M_3$ to reach
a halting configuration from a blank tape.

\begin{thm}\label{Thm5.1}
We have the following transitions between configurations of Turing machine $M_3$.
Let

$\begin{array}{lcl}
C(n,1) & = & ^\omega 013^n(B0)0^\omega,\\
C(n,2) & = & ^\omega 023^n(B0)0^\omega,\\
C(n,3) & = & ^\omega 03^n(B0)0^\omega,\\
C(n,4) & = & ^\omega 04113^n(B0)0^\omega,\\
C(n,5) & = & ^\omega 04123^n(B0)0^\omega,\\
C(n,6) & = & ^\omega 0413^n(B0)0^\omega,\\
C(n,7) & = & ^\omega 0423^n(B0)0^\omega,\\
C(n,8) & = & ^\omega 043^n(B0)0^\omega.
\end{array}$

\noindent Then

$\begin{array}{clll}
(a) & ^\omega 0(A0)0^\omega & \vdash(1) &  C(0,1),
\end{array}$

and, for all $k \ge 0$,

$\begin{array}{clll}
(b) &  C(2k,1)    & \vdash(3k^2 + 8k + 4)   & C(3k + 1,1),\\
(c) &  C(2k,2)    & \vdash(3k^2 + 14k + 9)  & C(3k + 2,1),\\
(d) &  C(2k,3)    & \vdash(3k^2 + 8k + 2)   & C(3k,1),\\
(e) &  C(2k,4)    & \vdash(3k^2 + 8k + 8)   & C(3k + 3,1),\\
(f) &  C(2k,5)    & \vdash(3k^2 + 14k + 13) & C(3k + 4,1),\\
(g) &  C(2k,6)    & \vdash(3k^2 + 8k + 6)   & C(3k + 2,1),\\
(h) &  C(2k,7)    & \vdash(3k^2 + 14k + 11) & C(3k + 3,1),\\
(i) &  C(2k,8)    & \vdash(3k^2 + 8k + 4)   & C(3k + 1,1),\\
(j) &  C(2k + 1,1)& \vdash(3k^2 + 8k + 4)   & C(3k + 1,2),\\
(k) &  C(2k + 1,2)& \vdash(3k^2 + 8k + 4)   & C(3k + 2,3),\\
(l) &  C(2k + 1,3)& \vdash(3k^2 + 8k + 22)  & C(3k + 1,4),\\
(m) &  C(2k + 1,4)& \vdash(3k^2 + 8k + 4)   & C(3k + 1,5),\\
(n) &  C(2k + 1,5)& \vdash(3k^2 + 8k + 4)   & C(3k + 2,6),\\
(o) &  C(2k + 1,6)& \vdash(3k^2 + 8k + 4)   & C(3k + 1,7),\\
(p) &  C(2k + 1,7)& \vdash(3k^2 + 8k + 4)   & C(3k + 2,8),\\
(q) &  C(2k + 1,8)& \vdash(3k^2 + 5k + 3)   & ^\omega 01(H2)2^{3k}0^\omega.
\end{array}$
\end{thm}

\proof
A direct inspection of the transition table gives

$\begin{array}{clll}
(1)  & (A0)0    & \vdash(1)  & 1(B0),\\
(2)  & 0(A0)0   & \vdash(17) &  41(A0),\\
(3)  & (A0)2    & \vdash(1)  &  1(B2),\\
(4)  & 0(A1)    & \vdash(1)  &  (A0)2,\\
(5)  & 1(A1)    & \vdash(1)  &  (A1)2,\\
(6)  & 4(A1)    & \vdash(1)  &  (A4)2,\\
(7)  & (A2)0^2  & \vdash(2)  &  1^2(B0),\\
(8)  & (A2)2    & \vdash(1)  &  1(A2),\\
(9)  & 1(B0)    & \vdash(1)  &  (A1)0,\\
(10) & 3^2(B0)0 & \vdash(5)  &  41^2(B0),\\
(11) & (B2)0    & \vdash(1)  &  3(B0),\\
(12) & (B2)2    & \vdash(1)  &  3(B2),\\
(13) & 0(A4)    & \vdash(1)  &  (A0)2,\\
(14) & 1(A4)    & \vdash(1)  &  (A1)2,\\
(15) & 2(A4)    & \vdash(1)  &  (A2)2,\\
(16) & 0^23(A4) & \vdash(3)  &  (A0)02^2,\\
(17) & 13(A4)   & \vdash(4)  &  23(B2),\\
(18) & 23(A4)   & \vdash(4)  &  33(B2),\\
(19) & 3^2(A4)  & \vdash(4)  &  41(A2),\\
(20) & 43(A4)   & \vdash(3)  &  1(H2)2,\\
(21) & 04(A4)   & \vdash(2)  &  (A0)2^2.\\
\end{array}$

\noindent Iterating, respectively, (5), (8) and (12) gives

$\begin{array}{clll}
(22) & 1^k(A1) & \vdash(k)  & (A1)2^k,\\
(23) & (A2)2^k & \vdash(k) & 1^k(A2),\\
(24) & (B2)2^k & \vdash(k)  & 3^k(B2).\\
\end{array}$

\smallskip

\noindent Using (9), (22), (6), (19), (23) and (7) we get

$\begin{array}{clll}
 & 3^241^{k+1}(B0)0 & \vdash(2k + 9) & 41^{k+4}(B0),
\end{array}$

\noindent and the result is still true for $k = -1$, so we have

$\begin{array}{clll}
(25) & 3^241^k(B0)0 & \vdash(2k + 7) & 41^{k+3}(B0).
\end{array}$

\smallskip

\noindent For any $n \ge 0$, by induction on $k$, using (25), we get

$\begin{array}{clll}
(26) & 3^{2k}41^n(B0)0^k & \vdash(3k^2+(2n+4)k) & 41^{3k+n}(B0),
\end{array}$

\noindent so we have, for $n = 2$ in (26),

$\begin{array}{clll}
(27) & 3^{2k}41^2(B0)0^k & \vdash(3k^2+8k) & 41^{3k+2}(B0).
\end{array}$

\smallskip

\noindent Using (10), (27), (9), (22) and (6), we get

$\begin{array}{clll}
(28) & 3^{2k+2}(B0)0^{k+1} & \vdash(3k^2 + 11k + 8) & (A4)2^{3k+2}0.
\end{array}$

\smallskip

\noindent Using (3), (24) and (11) we get

$\begin{array}{clll}
 & (A0)2^{k+1}0 & \vdash(k + 2) & 13^{k+1}(B0),
\end{array}$

\noindent and the result is still true for $k = -1$, so we have

$\begin{array}{clll}
(29) & (A0)2^k0 & \vdash(k + 1) & 13^k(B0).
\end{array}$

\smallskip

\noindent Using (15), (23), (7), (9) and (22), we get

$\begin{array}{clll}
(30) & 2(A4)2^k0^2 & \vdash(2k + 7) & (A1)2^{k+2}0.
\end{array}$

\smallskip

We are now ready to prove the theorem.

\smallskip

\noindent Using (28), (14), (4) and (29) we get

$\begin{array}{clll}
 & 013^{2k+2}(B0)0^{k+1} & \vdash(3k^2 + 14k + 15) & 13^{3k+4}(B0),
\end{array}$

\noindent and the result is still true for $k = -1$, so we have

$\begin{array}{clll}
(31) & 013^{2k}(B0)0^k & \vdash(3k^2 + 8k + 4) & 13^{3k+1}(B0).
\end{array}$

\smallskip

\noindent Using (28), (30), (4) and (29) we get

$\begin{array}{clll}
 & 023^{2k+2}(B0)0^{k+2} & \vdash(3k^2 + 20k + 26) & 13^{3k+5}(B0),
\end{array}$

\noindent and the result is still true for $k = -1$, so we have

$\begin{array}{clll}
(32) & 023^{2k}(B0)0^{k+1} & \vdash(3k^2 + 14k + 9) & 13^{3k+2}(B0).
\end{array}$

\smallskip

\noindent Using (28), (13) and (29) we get

$\begin{array}{clll}
 & 03^{2k+2}(B0)0^{k+1} & \vdash(3k^2 + 14k + 13) & 13^{3k+3}(B0),
\end{array}$

\noindent and the result is still true for $k = -1$, so we have

$\begin{array}{clll}
(33) & 03^{2k}(B0)0^k & \vdash(3k^2 + 8k + 2) & 13^{3k}(B0).
\end{array}$

\smallskip

\noindent Using (28), (14), (5), (6), (13) and (29) we get

$\begin{array}{clll}
 & 04113^{2k+2}(B0)0^{k+1} & \vdash(3k^2 + 14k + 19) & 13^{3k+6}(B0),
\end{array}$

\noindent and the result is still true for $k = -1$, so we have

$\begin{array}{clll}
(34) & 04113^{2k}(B0)0^k & \vdash(3k^2 + 8k + 8) & 13^{3k+3}(B0).
\end{array}$

\smallskip

\noindent Using (28), (30), (5), (6), (13) and (29) we get

$\begin{array}{clll}
 & 04123^{2k+2}(B0)0^{k+2} & \vdash(3k^2 + 20k + 30) & 13^{3k+7}(B0),
\end{array}$

\noindent and the result is still true for $k = -1$, so we have

$\begin{array}{clll}
(35) & 04123^{2k}(B0)0^{k+1} & \vdash(3k^2 + 14k + 13) & 13^{3k+4}(B0).
\end{array}$

\smallskip

\noindent Using (28), (14), (6), (13) and (29) we get

$\begin{array}{clll}
 & 0413^{2k+2}(B0)0^{k+1} & \vdash(3k^2 + 14k + 17) & 13^{3k+5}(B0),
\end{array}$

\noindent and the result is still true for $k = -1$, so we have

$\begin{array}{clll}
(36) & 0413^{2k}(B0)0^k & \vdash(3k^2 + 8k + 6) & 13^{3k+2}(B0).
\end{array}$

\smallskip

\noindent Using (28), (30), (6), (13) and (29) we get

$\begin{array}{clll}
 & 0423^{2k+2}(B0)0^{k+2} & \vdash(3k^2 + 20k + 28) & 13^{3k+6}(B0),
\end{array}$

\noindent and the result is still true for $k = -1$, so we have

$\begin{array}{clll}
(37) & 0423^{2k}(B0)0^{k+1} & \vdash(3k^2 + 14k + 11) & 13^{3k+3}(B0).
\end{array}$

\smallskip

\noindent Using (28), (21) and (29) we get

$\begin{array}{clll}
 & 043^{2k+2}(B0)0^{k+1} & \vdash(3k^2 + 14k + 15) & 13^{3k+4}(B0),
\end{array}$

\noindent and the result is still true for $k = -1$, so we have

$\begin{array}{clll}
(38) & 043^{2k}(B0)0^k & \vdash(3k^2 + 8k + 4) & 13^{3k+1}(B0).
\end{array}$

\smallskip

\noindent Using (28), (17), (24) and (11) we get

$\begin{array}{clll}
 & 13^{2k+3}(B0)0^{k+1} & \vdash(3k^2 + 14k + 15) & 23^{3k+4}(B0),
\end{array}$

\noindent and the result is still true for $k = -1$, so we have

$\begin{array}{clll}
(39) & 13^{2k+1}(B0)0^k & \vdash(3k^2 + 8k + 4) & 23^{3k+1}(B0).
\end{array}$

\smallskip

\noindent Using (28), (18), (24) and (11) we get

$\begin{array}{clll}
 & 23^{2k+3}(B0)0^{k+1} & \vdash(3k^2 + 14k + 15) & 3^{3k+5}(B0),
\end{array}$

\noindent and the result is still true for $k = -1$, so we have

$\begin{array}{clll}
(40) & 23^{2k+1}(B0)0^k & \vdash(3k^2 + 8k + 4) & 3^{3k+2}(B0).
\end{array}$

\smallskip

\noindent Using (28), (16), (2) and (29) we get

$\begin{array}{clll}
 & 0^33^{2k+3}(B0)0^{k+1} & \vdash(3k^2 + 14k + 33) & 41^23^{3k+4}(B0),
\end{array}$

\noindent and the result is still true for $k = -1$, so we have

$\begin{array}{clll}
(41) & 0^33^{2k+1}(B0)0^k & \vdash(3k^2 + 8k + 22) & 41^23^{3k+1}(B0).
\end{array}$

\smallskip

\noindent Using (39) we get

$\begin{array}{clll}
(42) & 04113^{2k+1}(B0)0^k & \vdash(3k^2 + 8k + 4) & 04123^{3k+1}(B0).
\end{array}$

\smallskip

\noindent Using (40) we get

$\begin{array}{clll}
(43) & 04123^{2k+1}(B0)0^k & \vdash(3k^2 + 8k + 4) & 0413^{3k+2}(B0).
\end{array}$

\smallskip

\noindent Using (39) we get

$\begin{array}{clll}
(44) & 0413^{2k+1}(B0)0^k & \vdash(3k^2 + 8k + 4) & 0423^{3k+1}(B0).
\end{array}$

\smallskip

\noindent Using (40) we get

$\begin{array}{clll}
(45) & 0423^{2k+1}(B0)0^k & \vdash(3k^2 + 8k + 4) & 043^{3k+2}(B0).
\end{array}$

\smallskip

\noindent Using (28) and (20) we get

$\begin{array}{clll}
 & 43^{2k+3}(B0)0^{k+1} & \vdash(3k^2 + 11k + 11) & 1(H2)2^{3k+3}0,
\end{array}$

\noindent and the result is still true for $k = -1$, so we have

$\begin{array}{clll}
(46) & 43^{2k+1}(B0)0^k & \vdash(3k^2 + 5k + 3) & 1(H2)2^{3k}0.
\end{array}$

\smallskip

\noindent Results (1) and (31)--(46) give results (a)--(p) of the theorem.
\qed

\smallskip

\noindent Using the rules of this theorem, we have, in 2002 transitions,
$$^\omega0(A0)0^\omega \vdash(1)\quad C(0,1) \vdash(4)\quad C(1,1) \vdash(\ ) \cdots
\vdash(\ )\quad \mbox{END}.$$

\noindent Let $g_3$ be the pure Collatz-like function with parameter defined by: for $k \ge 0$,

$$\begin{array}{rcl|rcl}
g_3(2k, 1) & = & (3k + 1, 1) \quad&\quad g_3(2k + 1, 1) & = & (3k + 1, 2)\\
g_3(2k, 2) & = & (3k + 2, 1) & g_3(2k + 1, 2) & = & (3k + 2, 3)\\
g_3(2k, 3) & = & (3k, 1) & g_3(2k + 1, 3) & = & (3k + 1, 4)\\
g_3(2k, 4) & = & (3k + 3, 1) & g_3(2k + 1, 4) & = & (3k + 1, 5)\\
g_3(2k, 5) & = & (3k + 4, 1) & g_3(2k + 1, 5) & = & (3k + 2, 6)\\
g_3(2k, 6) & = & (3k + 2, 1) & g_3(2k + 1, 6) & = & (3k + 1, 7)\\
g_3(2k, 7) & = & (3k + 3, 1) & g_3(2k + 1, 7) & = & (3k + 2, 8)\\
g_3(2k, 8) & = & (3k + 1, 1) & g_3(2k + 1, 8) &   & \mbox{undefined}
\end{array}$$
Then $g_3^{2001}(0,1)$ is undefined.

\begin{prop}\label{Prop5.2}
The behavior of Turing machine $M_3$, on inputs $02^n$, $n \ge 1$,
depends on the behavior of iterated $g_3^k(n,1)$, $k \ge 1$.
\end{prop}

\proof
We have $^\omega0(A0)2^n0^\omega \vdash(n+1)\quad
^\omega013^n(B0)0^\omega = C(n,1)$.\qed

\smallskip

Since the behavior of iterated $g_3^k(n,1)$ is an open problem in mathematics,
this si also the case for the halting problem for Turing machine $M_3$.

\smallskip

Note that the way by which a high score is obtained is particularly clear
for machine $M_3$. The parameter $p$, $1 \le p \le 8$, can be seen as a state.
If $n$ is odd, $g_3(n,p) = (n',p+1)$, the state goes from $p$ to $p+1$,
and the computation stops when state 8 is reached. If $n$ is even,
$g_3(n,p) = (n',1)$, and the state goes back to 1.

\section{Exponential Collatz-like}\label{Sect6}

Let $M_4$ be the $6 \times 2$ Turing machine defined by Table~\ref{m4}.

\begin{table}[t]\centering
\begin{tabular}{|c|c|c|}
\hline
$M_4$ &  0  &  1\\
\hline
$A$   & 1RB & 1LE\\
\hline
$B$   & 1RC & 1RF\\
\hline
$C$   & 1LD & 0RB\\
\hline
$D$   & 1RE & 0LC\\
\hline
$E$   & 1LA & 0RD\\
\hline
$F$   & 1RH & 1RC\\
\hline
\end{tabular}
\caption{\label{m4}Machine $M_4$ discovered in June 2010 by P.\ Kropitz}
\end{table}

We have $s(M_4) > 7.4 \times 10^{36534}$ and $\sigma(M_4) > 3.5 \times 10^{18267}$.

This machine is the current champion for the busy beaver competition for $6 \times 2$
machines. It was discovered in June 2010 by Pavel Kropitz.

The following theorem gives the rules observed by Turing machine $M_4$.

\begin{thm}\label{Thm6.1}
We have the following transitions between configurations of Turing machine $M_4$.
Let $C(n)$ = $^\omega 0(A0)1^n0^\omega$. Then

$\begin{array}{clll}
(a) & C(0) & \vdash(29) &  C(9),\\
(b) & C(2) & \vdash(36) &  C(11),\\
(c) & C(3) & \vdash(48) &  C(13),
\end{array}$

and, for all $k \ge 0$,

$\begin{array}{clll}
(d) &  C(3k + 1)  & \vdash(3k + 3)   & ^\omega0111(011)^k(H0)0^\omega,\\
(e) &  C(9k + 5)  & \vdash((4802\times 16^{k+1} + 6370\times 4^{k+1} + 2280k - 25362)/270)  & C((98\times 4^k - 11)/3),\\
(f) &  C(9k + 6)  & \vdash((125\times 16^{k+2} - 575\times 4^{k+2} + 228k - 2226)/27)  & C((50\times 4^{k+1} - 59)/3),\\
(g) &  C(9k + 8)  & \vdash((4802\times 16^{k+1} + 6370\times 4^{k+1} + 2280k - 11592)/270)  & C((98\times 4^k + 1)/3),\\
(h) &  C(9k + 9)  & \vdash((125\times 16^{k+2} + 325\times 4^{k+2} + 228k - 2289)/27)  & C((50\times 4^{k+1} - 11)/3),\\
(i) &  C(9k + 11) & \vdash((4802\times 16^{k+2} - 11270\times 4^{k+2} + 2280k - 22452)/270)  & C((98\times 4^{k+1} - 59)/3),\\
(j) &  C(9k + 12)  & \vdash((125\times 16^{k+2} + 325\times 4^{k+2} + 228k - 912)/27)  & C((50\times 4^{k+1} + 1)/3).
\end{array}$

Note that the behavior of this Turing machine on the blank tape involves only items \emph{(a), (d), (h)} and \emph{(j)}.
\end{thm}

\proof
A direct inspection of the transition table gives

$\begin{array}{clll}
(1)  & 0^3(A0)0^6    & \vdash(29)  & (A0)1^9,\\
(2)  & 0^4(A0)1^20^5 & \vdash(36)  & (A0)1^{11},\\
(3)  & 0^4(A0)1^30^6 & \vdash(48)  & (A0)1^{13},\\
(4)  & (A0)1        & \vdash(1)   & 1(B1),\\
(5)  & 01(E0)       & \vdash(2)   & (E0)11,\\
(6)  & 0(E0)        & \vdash(1)   & (A0)1,\\
(7)  & 11(E0)       & \vdash(2)   & (E1)11,\\
(8)  & 01(C0)       & \vdash(2)   & (C0)01,\\
(9)  & 11(C0)       & \vdash(2)   & (C1)01,\\
(10) & 0(C0)        & \vdash(2)   & 1(E1),\\
(11) & (E1)01       & \vdash(2)   & 01(E1),\\
(12) & (E1)00       & \vdash(2)   & 01(E0),\\
(13) & (E1)1        & \vdash(2)   & (C0)0,\\
(14) & (B1)1^3      & \vdash(3)   & 110(B1),\\
(15) & (B1)00       & \vdash(2)   & 11(H0),\\
(16) & (B1)10       & \vdash(6)   & 01(C1),\\
(17) & (B1)1100     & \vdash(12)  & (01)^2(C1),\\
(18) & (C1)01       & \vdash(2)   & 01(C1),\\
(19) & (C1)00       & \vdash(2)   & 01(C0),\\
(20) & (C1)1^60^6   & \vdash(44)  & 01(E1)1^{10},\\
(21) & (C1)1^80^{11} & \vdash(113) & 1(01)^5(E1)1^8.
\end{array}$

\smallskip

\noindent Iterating, respectively, (5), (8), (11), (14) and (18) gives

$\begin{array}{clll}
(22) & (01)^k(E0) & \vdash(2k)  & (E0)1^{2k},\\
(23) & (01)^k(C0) & \vdash(2k)  & (C0)(01)^k,\\
(24) & (E1)(01)^k & \vdash(2k)  & (01)^k(E1),\\
(25) & (B1)1^{3k}  & \vdash(3k)  & (110)^k(B1),\\
(26) & (C1)(01)^k & \vdash(2k)  & (01)^k(C1).
\end{array}$

\smallskip

\noindent Using (19), (23), (10) and (24), we get

$\begin{array}{clll}
(27) & 0(01)^k(C1)00 & \vdash(4k + 8) & 1(01)^{k+1}(E1).
\end{array}$

\smallskip

\noindent Using (12), (22) and (6), we get

$\begin{array}{clll}
(28) & 0(01)^k(E1)00 & \vdash(2k + 5) & (A0)1^{2k+3}.
\end{array}$

\smallskip

\noindent Using (12), (22) and (7), we get

$\begin{array}{clll}
(29) & 11(01)^k(E1)00 & \vdash(2k + 6) & (E1)1^{2k+4}.
\end{array}$

\smallskip

\noindent Using (13), (23), (10) and (24), we get

$\begin{array}{clll}
(30) & 0(01)^k(E1)1 & \vdash(4k + 4) & 1(01)^k(E1)0.
\end{array}$

\smallskip

\noindent Using (13), (23), (9) and (26), we get

$\begin{array}{clll}
(31) & 11(01)^k(E1)1 & \vdash(4k+6) & (01)^{k+1}(C1)0.
\end{array}$

\smallskip

\noindent Using (20), (30), (11), (31) and (18), we get

$\begin{array}{clll}
(32) & 10(01)^k(C1)1^60^6 & \vdash(8k+70) & (01)^{k+4}(C1)1^6.
\end{array}$

\smallskip

\noindent By induction on $n$, using (32),  we get

$\begin{array}{clll}
(33) & (10)^n(01)^k(C1)1^60^{6n} & \vdash(16n^2 + 8kn + 54n) & (01)^{4n+k}(C1)1^6.
\end{array}$

\smallskip

\noindent Using (30) and (11), we get

$\begin{array}{clll}
(34) & 00(01)^k(E1)11 & \vdash(4k + 6) & (01)^{k+2}(E1).
\end{array}$

\smallskip

\noindent By induction on $n$, using (34), we get

$\begin{array}{clll}
(35) & 0^{2n}(01)^k(E1)1^{2n} & \vdash(4n^2 + 4kn + 2n) & (01)^{2n+k}(E1).
\end{array}$

\smallskip

\noindent Using (20), (31), (18), (21), (31), (18) and (33), we get

$\begin{array}{clll}
(36) & 11(01)^k(C1)1^60^{6k+29} & \vdash(16k^2 + 178k + 481) & 0(01)^{4k+15}(C1)1^6.
\end{array}$

\smallskip

\noindent Using (20), (35) and (28), we get

$\begin{array}{clll}
(37) & 0^{11}(01)^k(C1)1^60^8 & \vdash(22k + 201) & (A0)1^{2k+25}.
\end{array}$

\smallskip

\noindent Using (20), (35), (31), (18), (32) and (36), we get

$\begin{array}{clll}
(38) & (110)^30(01)^k(C1)1^60^{6k+101} & \vdash(16k^2 + 514k + 4045) & 0(01)^{4k+55}(C1)1^6.
\end{array}$

\smallskip

\noindent By induction on $k$, using (38) we get

$\begin{array}{clll}
(39) & (110)^{3k}0(01)^n(C1)1^60^a & \vdash(T) & 0(01)^b(C1)1^6,
\end{array}$

\noindent with $a = (2(3n+55)4^k-6n-27k-110)/3$, $b = ((3n+55)4^k-55)/3$, and
$T=\frac{16(3n+55)^2}{135}16^k-\frac{218(3n+55)}{27}4^k-\frac{5}{9}k-\frac{16(3n+55)^2}{135}+\frac{218(3n+55)}{27}$.

\smallskip

\noindent Using (4), (25), (17), (27), (29), (31), (18), (21), (31), (18) and (33), we get

$\begin{array}{clll}
(40) & (A0)1^{3k+9}0^{29} & \vdash(3k +484) & 1(110)^k0(01)^{15}(C1)1^6.
\end{array}$

\smallskip

\noindent Using (40), (39), (32) and (37), we get

$\begin{array}{clll}
(41) & 0^{11}(A0)1^{9k+9}0^a & \vdash(T) & (A0)1^b,
\end{array}$

\noindent with $a = (50 \times 4^{k+1} - 11)/3 - 9k - 20$, $b = (50 \times 4^{k+1} - 11)/3$, and
$T = (125 \times 16^{k+2} + 325 \times 4^{k+2} + 228k - 2289)/27$.

\smallskip

\noindent Using (40), (39), (20), (35), (31), (18) and (37), we get

$\begin{array}{clll}
(42) & 0^{12}(A0)1^{9k+12}0^a & \vdash(T) & (A0)1^b,
\end{array}$

\noindent with $a = (50 \times 4^{k+1} + 1)/3 - 9k - 24$, $b = (50 \times 4^{k+1} + 1)/3$, and
$T = (125 \times 16^{k+2} + 325 \times 4^{k+2} + 228k - 912)/27$.

\smallskip

\noindent Using (4), (25) and (15), we get

$\begin{array}{clll}
(43) & (A0)1^{3k+1}00 & \vdash(3k + 3) & 111(011)^k(H0).
\end{array}$

\smallskip

Results (1), (43), (41) and (42) are results (a), (d), (h) and (j) of the theorem.
They are sufficient to analyze the behavior of the Turing machine on a blank tape.
The following gives its behavior from configurations $^\omega0(A0)1^n0^\omega$, where $n = 9k + m$,
$m \in \{5,6,8,11\}$.

\smallskip

\noindent Using (4), (25), (16), (27), (29), (31), (18), (32) and (36), we get

$\begin{array}{clll}
(44) & (A0)1^{3k+14}0^{88} & \vdash(3k + 3076) & 1(110)^k0(01)^{47}(C1)1^6.
\end{array}$

\smallskip

\noindent Using (44), (39), (32) and (37), we get

$\begin{array}{clll}
 & 0^{11}(A0)1^{9k+14}0^a & \vdash(T) & (A0)1^b,
\end{array}$

\noindent with $b = (98 \times 4^{k+1} - 11)/3$ and
$T = (4802 \times 16^{k+2} + 6370 \times 4^{k+2} + 2280(k+1) - 25362)/270$,

\noindent and this result is still true for $k = -1$, since

$\begin{array}{clll}
 & 0^{11}(A0)1^5 0^{13} & \vdash(285) & (A0)1^{29},
\end{array}$

\noindent so we get

$\begin{array}{clll}
(45) & 0^{11}(A0)1^{9k+5}0^a & \vdash(T) & (A0)1^b,
\end{array}$

\noindent with $b = (98 \times 4^k - 11)/3$ and
$T = (4802 \times 16^{k+1} + 6370 \times 4^{k+1} + 2280k - 25362)/270$.

\smallskip

\noindent Using (20), (35), (31), (18), (32), (36) and (37), we get

$\begin{array}{clll}
(46) & 0^{10}1(110)^20(01)^k(C1)1^60^{6k+103} & \vdash(16k^2 + 570k + 4886) & (A0)1^{8k+127}.
\end{array}$

\smallskip

\noindent Using (40), (39) and (46), we get

$\begin{array}{clll}
 & 0^{10}(A0)1^{9k+15}0^a & \vdash(T) & (A0)1^b,
\end{array}$

\noindent with $b = (50 \times 4^{k+2} - 59)/3$ and
$T = (125 \times 16^{k+3} - 575 \times 4^{k+3} + 228(k+1) - 2226)/27$,

\noindent and this result is still true for $k = -1$, since

$\begin{array}{clll}
 & 0^{10}(A0)1^6 0^{31} & \vdash(762) & (A0)1^{47},
\end{array}$

\noindent so we get

$\begin{array}{clll}
(47) & 0^{10}(A0)1^{9k+6}0^a & \vdash(T) & (A0)1^b,
\end{array}$

\noindent with $b = (50 \times 4^{k+1} - 59)/3$ and
$T = (125 \times 16^{k+2} - 575 \times 4^{k+2} + 228k - 2226)/27$.

\smallskip

\noindent Using (20), (30), (11), (31), (18) and (37), we get

$\begin{array}{clll}
(48) & 0^{12}11100(01)^k(C1)1^60^{14} & \vdash(30k + 407) & (A0)1^{2k+37}.
\end{array}$

\smallskip

\noindent Using (44), (39) and (48), we get

$\begin{array}{clll}
 & 0^{12}(A0)1^{9k+17}0^a & \vdash(T) & (A0)1^b,
\end{array}$

\noindent with $b = (98 \times 4^{k+1} + 1)/3$ and
$T = (4802 \times 16^{k+2} + 6370 \times 4^{k+2} +2280(k+1) - 11592)/270$,

\noindent and this result is still true for $k = -1$, since

$\begin{array}{clll}
 & 0^{12}(A0)1^8 0^{13} & \vdash(336) & (A0)1^{33},
\end{array}$

\noindent so we get

$\begin{array}{clll}
(49) & 0^{12}(A0)1^{9k+8}0^a & \vdash(T) & (A0)1^b,
\end{array}$

\noindent with $b = (98 \times 4^k + 1)/3$ and
$T = (4802 \times 16^{k+1} + 6370 \times 4^{k+1} +2280k - 11592)/270$.

\smallskip

\noindent Using (44), (39) and (46), we get

$\begin{array}{clll}
 & 0^{10}(A0)1^{9k+20}0^a & \vdash(T) & (A0)1^b,
\end{array}$

\noindent with $b = (98 \times 4^{k+2} - 59)/3$ and
$T = (4802 \times 16^{k+3} - 11270 \times 4^{k+3} +2280(k+1) - 22452)/270$,

\noindent and this result is still true for $k = -1$, since

$\begin{array}{clll}
 & 0^{10}(A0)1^{11} 0^{90} & \vdash(3802) & (A0)1^{111},
\end{array}$

\noindent so we get

$\begin{array}{clll}
(50) & 0^{10}(A0)1^{9k+11}0^a & \vdash(T) & (A0)1^b,
\end{array}$

\noindent with $b = (98 \times 4^{k+1} - 59)/3$ and
$T = (4802 \times 16^{k+2} - 11270 \times 4^{k+2} +2280k - 22452)/270$.

\smallskip

Results (45), (47), (49) and (50) are results (e), (f), (g) and (i)
of the theorem.

Results (1), (2) and (3) give special cases (a), (b) and (c) of the theorem. 
\qed

\smallskip

\noindent Using the rules of this theorem, we have,
$$^\omega0(A0)0^\omega \vdash(29)\quad C(9) \vdash(1293)\quad C(63) \vdash(19,884,896,677 )$$
$$C(273063) \vdash(T_1)\quad C((50 \times 4^{30340} + 1)/3) \vdash(T_2)\quad
^\omega0111(011)^K(H0)0^\omega,$$
with $T_1 = (125 \times 16^{30341} + 325 \times 4^{30341} + 6916380)/27$,
$T_2 = (50 \times 4^{30340} + 7)/3$, $K = (50 \times 4^{30340} - 2)/9$.

The total time is $s(M_4) = (125 \times 16^{30341} + 1750 \times 4^{30340} + 15)/27
+ 19,885,154,163$, and the final number of symbols 1 is $\sigma(M_4) = 
(25 \times 4^{30341} + 23)/9$.

Let us give some informal comments on how Turing machine $M_4$ works.
As in the case of machine $M_1$, edge effects give the multiplicity of cases.
Also as in the case of machine $M_1$, machine $M_4$ nibbles strings $(01)^n$,
transforming them into strings four times longer (see transitions (32) and
(33) in the proof of Theorem \ref{Thm6.1}). A new idea allows machine $M_4$
to reach exponential scores. Machine $M_4$ nibbles strings $(110)^{3k}$ and,
each time a string $(110)^3$ is nibbled, the length of a string $(01)^n$
becomes four times longer (see transitions (38) and (39) in the proof of Theorem \ref{Thm6.1}).
This explains why powers of four are obtained in the transitions of Theorem \ref{Thm6.1}.

\smallskip

Let $g_4$ be the exponential Collatz-like function defined by: for $k \ge 0$,

$$\begin{array}{rcl}
g_4(0) & = & 9,\\
g_4(2) & = & 11,\\
g_4(3) & = & 13,\\
g_4(3k + 1) &   & \mbox{undefined},\\
g_4(9k + 5) & = & (98 \times 4^k - 11)/3,\\
g_4(9k + 6) & = & (50 \times 4^{k+1} - 59)/3,\\
g_4(9k + 8) & = & (98 \times 4^k + 1)/3,\\
g_4(9k + 9) & = & (50 \times 4^{k+1} - 11)/3,\\
g_4(9k + 11) & = & (98 \times 4^{k+1} - 59)/3,\\
g_4(9k + 12) & = & (50 \times 4^{k+1} + 1)/3.
\end{array}$$
Then $g_4^{5}(0)$ is undefined.

The theorem gives immediately the following proposition.

\begin{prop}\label{Prop6.2}
The behavior of Turing machine $M_4$, on inputs $01^n$, $n \ge 0$,
depends on the behavior of iterated $g_4^k(n)$, $k \ge 1$.\qed
\end{prop}

\noindent Since the behavior of iterated $g_4^k(n)$ is an open problem in mathematics,
this is also the case for the halting problem for Turing machine $M_4$.

Let $h_4(n) = \min\{k:g_4^k(n) \mbox{\ is undefined}\}$. We have seen that
$h_4(0) = 5$. We also have $h_4(2) = 8$, and $C(2) \vdash(T)$ END with
$T> 10^{10^{10^{10^{18641000}}}}$. We also have $h_4(36) = 15$.

\section{Unclassifiable machine}\label{Sect7}

Let $M_5$ be the $6 \times 2$ Turing machine defined by Table~\ref{m5}.

\begin{table}[t]\centering
\begin{tabular}{|c|c|c|}
\hline
$M_5$ &  0  &  1\\
\hline
$A$   & 1RB & 0LD\\
\hline
$B$   & 1RC & 0RF\\
\hline
$C$   & 1LC & 1LA\\
\hline
$D$   & 0LE & 1RH\\
\hline
$E$   & 1LA & 0RB\\
\hline
$F$   & 0RC & 0RE\\
\hline
\end{tabular}
\caption{\label{m5}Machine $M_5$ discovered in May 2010 by P.\ Kropitz}
\end{table}

We have $s(M_5) > 3.8 \times 10^{21132}$ and $\sigma(M_5) > 3.1 \times 10^{10566}$.

This machine was discovered in May 2010 by Pavel Kropitz.
It was the champion for the busy beaver competition for $6 \times 2$
machines from May to June 2010.

The following theorem is adapted from an analysis of S. Ligocki \cite{Li3}.
It gives the rules that enable Turing machine $M_5$ to reach
a halting configuration from a blank tape.

\begin{thm}\label{Thm7.1}
We have the following transitions between configurations of Turing machine $M_5$.
Let $C(k,n)$ = $^\omega010^n1(C1)1^{3k}0^\omega$. Then

$\begin{array}{clll}
(a) & ^\omega0(A0)0^\omega & \vdash(47) & C(2,5),
\end{array}$

and, for all $k \ge 0$,

$\begin{array}{clll}
(b) & C(k,0)   & \vdash(3) & ^\omega01(H0)1^{3k+1}0^\omega,\\
(c) & C(k,1)   & \vdash(3k + 37) & C(2,3k+2),\\
(d) & C(k,2)   & \vdash(12k + 44) & C(k+2,4),\\
(e) & C(k,3)   & \vdash(3k + 57) & C(2,3k+8),\\
(f) & C(k,n+4) & \vdash(27k^2 + 105k + 112) & C(3k+5,n).\\
\end{array}$
\end{thm}

\proof
A direct inspection of the transition table gives

$\begin{array}{clll}
(1)  & 0^4(A0)0^9   & \vdash(47)  & 10^51(C1)1^6,\\
(2)  & 0(C0)        & \vdash (1)  & (C0)1,\\
(3)  & 1(C0)        & \vdash (1)  & (C1)1,\\
(4)  & (B1)00       & \vdash (2)  & 0^2(C0),\\
(5)  & (B1)01       & \vdash (4)  & 01(B1),\\
(6)  & (B1)10       & \vdash (4)  & 01(B1),\\
(7)  & (B1)10^3     & \vdash (14) & 10^3(B1),\\
(8)  & (B1)1^20^2   & \vdash (7)  & 0^21(B1)1,\\
(9)  & (B1)1^3      & \vdash (3)  & 0^3(B1),\\
(10) & 0(C1)        & \vdash (2)  & 1(B1),\\
(11) & 11(C1)       & \vdash (3)  & 1(H0)1,\\
(12) & 0^31(C1)     & \vdash (10) & (C1)1^4,\\
(13) & 0^310^21(C1) & \vdash (8)  & 1(B1)0^210^21,\\
(14) & 101(C1)      & \vdash (8)  & 1(B1)1^2.
\end{array}$

\smallskip

\noindent Iterating, respectively, (2) and (9) gives

$\begin{array}{clll}
(15) & 0^k(C0)   & \vdash (k)  & (C0)1^k,\\
(16) & (B1)1^{3k} & \vdash (3k) & 0^{3k}(B1).
\end{array}$

\smallskip

\noindent Using (4), (15) and (3), we get

$\begin{array}{clll}
(17) & 10^k(B1)00 & \vdash(k + 5) & (C1)1^{k+3}.
\end{array}$

\smallskip

\noindent Using (16), (17) and (10), we get

$\begin{array}{clll}
(18) & 01(B1)1^{3k}00 & \vdash(6k + 7) & 1(B1)1^{3k+3}.
\end{array}$

\smallskip

\noindent By induction on $k$, using (18), we get

$\begin{array}{clll}
(19) & 0^k1(B1)0^{2k} & \vdash(3k^2 + 4k) & 1(B1)1^{3k}.
\end{array}$

\smallskip

\noindent Using (12), (10), (16), (6), (19), (16) and (17), we get

$\begin{array}{clll}
(20) & 0^41(C1)1^{3k}0^{6k+11} & \vdash(27k^2 + 105k + 112) & 1(C1)1^{9k+15}.
\end{array}$

\smallskip

\noindent Using (14), (16), (8), (7) and (17), we get

$\begin{array}{clll}
(21) & 101(C1)1^{3k}0^7 & \vdash(3k + 37) & 10^{3k+2}1(C1)1^6.
\end{array}$

\smallskip

\noindent Using (13), (17), (10), (9), (6), (5), (16), (17), (10), (16) and (17), we get

$\begin{array}{clll}
(22) & 0^410^21(C1)1^{3k}0^4 & \vdash(12k + 44) & 10^41(C1)1^{3k+6}.
\end{array}$

\smallskip

\noindent Using (12), (12), (10), (16), (8), (6), (17), (10), (9) and (17), we get

$\begin{array}{clll}
(23) & 0^410^31(C1)1^{3k}0^7 & \vdash(3k + 57) & 10^{3k+8}1(C1)1^6.
\end{array}$

\smallskip

\noindent Results (1), (11), (21), (22), (23) and (20) give results (a)--(f) of the theorem.
\qed

\smallskip

\noindent Using the rules of this theorem, we have, in 22158 transitions,
$$^\omega0(A0)0^\omega \vdash(47)\quad C(2,5) \vdash(430)\quad C(11,1) \vdash(\ ) \cdots
\vdash(\ )\quad \mbox{END}.$$
Let us give some informal comments on how Turing machine $M_5$ works.
The computation of machine $M_5$ on a blank tape is a succession
of two types of phases. In the first one, the second parameter $p$ of
configuration $C(n,p)$ is smaller than four, and transitions (b)--(e)
of Theorem \ref{Thm7.1} immediately give a greater value to this second parameter.
In the second type of phase, the second parameter is greater than four,
and slowly decreases, while the first parameter rapidly increases.
In this second type of phase, machine $M_5$ nibbles strings of symbols 0
two by two, and each time a string of two symbols 0 is nibbled,
a string of three symbols 1 is produced (see transitions (18) and (19)
in the proof of Theorem \ref{Thm7.1}). In the computation of Turing machine $M_5$
on a blank tape, phases of the second type occur four times:
\begin{itemize}
\item[-] from $C(2,5)$ to $C(11,1)$ (1 transition),
\item[-] from $C(2,35)$ to $C(29522,3)$ (8 transitions),
\item[-] from $C(2,88574)$ to $C(m,2)$ (22143 transitions),
\item[-] from $C(m+2,4)$ to $C(3(m+2)+5,0)$ (1 transition),
\end{itemize}
where $m$ is a big number.

\smallskip

Let $g_5$ be the partial function defined by: for $k,n \ge 0$,

$$\begin{array}{rcl}
g_5(k, 0) &   & \mbox{undefined},\\
g_5(k, 1) & = & (2, 3k + 2),\\
g_5(k, 2) & = & (k + 2, 4),\\
g_5(k, 3) & = & (2, 3k + 8),\\
g_5(k, n + 4) & = & (3k + 5, n).
\end{array}$$
Then $g_5^{22157}(2,5)$ is undefined.

\begin{prop}\label{Prop7.2}
The behavior of Turing machine $M_5$, on inputs $01^{3n+3}$, $n \ge 0$,
depends on the behavior of iterated $g_5^k(2,3n + 2)$.
\end{prop}

\proof
We have $^\omega0(A0)1^{3n+3}0^\omega \vdash(3n + 30)\quad
^\omega01^{3n+2}1(C1)1^60^\omega = C(2,3n+2)$.\qed

\smallskip

Since the behavior of iterated $g_5^k(n,p)$ is an open problem in mathematics,
this is also the case for the halting problem for Turing machine $M_5$.

\smallskip

The following proposition shows that some configurations take a long time to halt.

\begin{prop}\label{Prop7.3}
For Turing machine $M_5$, we have
$C(9,1) \vdash(T)\ \ \mbox{END}$
with $T > 10^{10^{10^{10^{10^{3520}}}}}$.
\end{prop}

\proof
By induction on $n$, using Theorem 7.1 (f), it is easy to prove
that, if  $n \ge 0$, $0 \le r \le 3$, we have
$$C(2,4n + r) \vdash(t_n)\quad C(u_n,r),$$
with $u_n = (3^{n+2}-5)/2$ and $t_n = (3 \times 9^{n+3} - 80 \times 3^{n+3} +
584n - 27)/32$.

By induction on $k$, it is easy to prove that, if $k \ge 2$, we have
$$3^{2^{k-1}} \equiv 2^{k+1} + 1\ (\mbox{mod}\ 2^{k+2})$$
so the multiplicative order of 3 modulo $2^{k+2}$ is $2^k$
for $k \ge 1$. Thus we can prove that, for $k \ge 1$, $n,m \ge 0$,
we have 
$$n \equiv m \ (\mbox{mod}\ 2^k) \iff u_n \equiv u_m \ (\mbox{mod}\ 2^{k+1}).$$
Now, suppose that, for $a \in \{1,3\}$, $n, n' \ge 1$, $q, q' \ge 1$, $0 \le r,r' \le 3$,
we have
$$C(n,a) \vdash(3n + 27 + 10a)\quad C(2,3n + 3a -1) = C(2, 4q + r)
\vdash(t_q)\quad C(u_q, r),$$
and
$$C(n',a) \vdash(3n' + 27 + 10a)\quad C(2,3n' + 3a -1) = C(2, 4q' + r')
\vdash(t_{q'})\quad C(u_{q'}, r'),$$
and let $k \ge 2$ such that $n \equiv n' \ (\mbox{mod}\ 2^{k+1})$.
Then it is easy to prove that $r = r'$ and $u_q \equiv u_{q'} \ (\mbox{mod}\ 2^k)$.
So the behavior of configurations $C(n,a)$ is mirrored by the behavior
of configurations $C(n',a)$ with $n' \le 2^k$ for suitable $k$.

In the following computation on $C(9,1)$:
$$\begin{array}{llll}
C(9, 1)       & \vdash(\ ) & C(2,4 \times 7 + 1)      & \vdash(t_7)\\
C(9839, 1)    & \vdash(\ ) & C(2, 4 \times 7379 + 3)  & \vdash(t_{7379})\\
C(u_{7379}, 3)  & \vdash(\ ) & C(2, 4 \times q_3 + r_3) & \vdash(t_{q_3})\\
C(u_{q_3}, r_3) & \vdash(\ ) & C(2, 4 \times q_4 + r_4) & \vdash(t_{q_4})\\
C(u_{q_4}, r_4) & \vdash(\ ) & C(2, 4 \times q_5 + r_5) & \vdash(t_{q_5})\\
C(u_{q_5}, r_5) & \vdash(\ ) & C(2, 4 \times q_6 + r_6) & \vdash(t_{q_6})\\
C(u_{q_6}, r_6) & \vdash(3)  & \mbox{END}              & 
\end{array}$$
we know that $r_6 = 0$ because we have
$$\begin{array}{llll}
u_{q_1} = u_7 \equiv 47\ (\mbox{mod}\ 64), & (3 \times 47) + 2 = (4 \times 35) + 3, & q_2' = 35, & r_2 = 3,\\
u_{q_2} \equiv u_{q_2'} \equiv 23\ (\mbox{mod}\ 32), & (3 \times 23) + 8 = (4 \times 19) + 1, & q_3' = 19, & r_3 = 1,\\
u_{q_3} \equiv u_{q_3'} \equiv 7\ (\mbox{mod}\ 16), & (3 \times 7) + 2 = (4 \times 5) + 3, & q_4' = 5, & r_4 = 3,\\
u_{q_4} \equiv u_{q_4'} \equiv 3\ (\mbox{mod}\ 8), & (3 \times 3) + 8 = (4 \times 4) + 1, & q_5' = 4, & r_5 = 1,\\
u_{q_5} \equiv u_{q_5'} \equiv 2\ (\mbox{mod}\ 4), & (3 \times 2) + 2 = (4 \times 2) + 0, & q_6' = 2, & r_6 = 0.
\end{array}$$
It is easy to see that, if $a \in \{1,3\}$, $n \ge 0$, and
$$C(n,a) \vdash(3n + 27 + 10a)\quad C(2,3n + 3a -1) = C(2, 4q + r)
\vdash(t_q)\quad C(u_q, r),$$
then $q \ge (3n - 1)/4$ and $u_q > (3^{3/4})^n > 2^n$.

And we also have $n \ge 5 \Rightarrow t_n > 68 \times 9^n$,
so, if $C(9, 1) \vdash(T)$ END, we have
$$T > t_{q_6} > 9^{q_6} > 9^{3u_{q_5}/4} > 5^{u_{q_5}},$$
and $u_{q_5} >2^{u_{q_4}}$, $u_{q_4} >2^{u_{q_3}}$, $u_{q_3} >2^{u_{q_2}} = 2^{u_{7379}}$,
so $T > 5^{2^{2^{2^{u_{7379}}}}}$.

Using $u_{7379} > 10^{3521}$, and, for $x \ge 1$, $2^{10^x} > 10^{10^{x-.53}}$,
$2^{10^{10^x}} > 10^{10^{10^{x-.03}}}$, $2^{10^{10^{10^x}}} > 10^{10^{10^{10^{x-.03}}}}$ and
$5^{10^{10^{10^{10^x}}}} > 10^{10^{10^{10^{10^{x-.03}}}}}$, we are done.\qed

\section{A potentially infinite set of rules}\label{Sect8}

Let $M_6$ be the $6 \times 2$ Turing machine defined by Table~\ref{m6}.

\begin{table}[t]\centering
\begin{tabular}{|c|c|c|}
\hline
$M_6$ &  0  &  1\\
\hline
$A$   & 1RB & 0RF\\
\hline
$B$   & 0LB & 1LC\\
\hline
$C$   & 1LD & 0RC\\
\hline
$D$   & 1LE & 1RH\\
\hline
$E$   & 1LF & 0LD\\
\hline
$F$   & 1RA & 0LE\\
\hline
\end{tabular}
\caption{\label{m6}Machine $M_6$ discovered in November 2007 by T.\ and S.\ Ligocki}
\end{table}

We have $s(M_6) > 8.9 \times 10^{1762}$ and $\sigma(M_6) > 2.5 \times 10^{881}$.

This machine was discovered in November 2007 by Terry and Shawn Ligocki.
It was the champion for the busy beaver competition for $6 \times 2$
machines from November to December 2007.

The complete analysis of Turing machine $M_6$ seems to need an infinite
set of rules, but proving this assertion could be difficult.
Of course, only a finite subset of these rules are needed when the machine
is launched on a blank tape.

The following theorem gives the rules that enable Turing machine $M_6$ to reach
a halting configuration from a blank tape.

Recall that bin($p$) is the usual binary representation of number $p$, and
R($w_1 \ldots w_n$) = $w_n \ldots w_1$.

\begin{thm}\label{Thm8.1}
We have the following transitions between configurations of Turing machine $M_6$.
Let
$C(n,p)$ = $^\omega0(F0)(10)^n\mbox{R(bin($p$))}0^\omega$,
so that $C(k,4m + 1) = C(k + 1,m)$.
Then

$\begin{array}{clll}
(a) & ^\omega0(A0)0^\omega & \vdash(6) & C(0,15),
\end{array}$

and, for all $k,m \ge 0$,

$\begin{array}{clll}
(b) & C(k,4m+3)         & \vdash(4k + 6)           & C(k+2,m),\\
(c) & C(2k,4m)          & \vdash(30k^2 + 20k + 15) & C(5k+2,2m+1),\\
(d) & C(2k+1,4m)        & \vdash(30k^2 + 40k + 25) & C(5k+2,32m+20),\\
(e) & C(k,8m+2)         & \vdash(8k + 20)          & C(k+3,2m+1),\\
(f) & C(2k,16m+6)       & \vdash(30k^2 + 40k + 23) & C(5k+2,32m+20),\\
(g) & C(2k+1,16m+6)     & \vdash(30k^2 + 80k + 63) & C(5k+7,2m+1),\\
(h) & C(k,32m+14)       & \vdash(4k + 18)          & C(k+3,2m+1),\\
(i) & C(2k,128m+94)     & \vdash(30k^2 + 40k + 39) & C(5k+2,256m+84),\\
(j) & C(2k+1,128m+94)   & \vdash(30k^2 + 80k + 79) & C(5k+9,m),\\
(k) & C(k,256m+190)     & \vdash(4k + 34)          & C(k+5,m),\\
(l) & C(k,512m+30)      & \vdash(2k + 43)          & ^\omega0(10)^k1(H0)(10)^2(01)^2\mbox{R(bin($m$))}0^\omega.
\end{array}$
\end{thm}

\proof
A direct inspection of the transition table gives

$\begin{array}{clll}
(1)  & 0^4(A0)0 & \vdash(6) & (F0)1^40,\\
(2)  & 0(F0)00  & \vdash(9) & (F0)1^3,\\
(3)  & (F0)10   & \vdash(2) & 10(F0),\\
(4)  & (F0)11   & \vdash(4) & (F1)10,\\
(5)  & (F0)01   & \vdash(4) & 10(C1),\\
(6)  & 10(E1)   & \vdash(2) & (E1)10,\\
(7)  & 1(E1)    & \vdash(2) & 1(H0),\\
(8)  & 0^3(E1)  & \vdash(3) & (F0)110,\\
(9)  & 00(F1)   & \vdash(2) & (F0)10,\\
(10) & 10(F1)   & \vdash(2) & (F1)10,\\
(11) & 1(F1)    & \vdash(1) & (E1)0,\\
(12) & (C1)1    & \vdash(1) & 0(C1),\\
(13) & (C1)0    & \vdash(1) & 0(C0),\\
(14) & 0^3(C0)  & \vdash(3) & (F0)1^3,\\
(15) & 100(C0)  & \vdash(3) & (F1)1^3.
\end{array}$

\smallskip

\noindent Iterating, respectively, (3), (10), (6) and (12) gives

$\begin{array}{clll}
(16) & (F0)(10)^k & \vdash(2k)  & (10)^k(F0),\\
(17) & (10)^k(F1) & \vdash(2k)  & (F1)(10)^k,\\
(18) & (10)^k(E1) & \vdash(2k)  & (E1)(10)^k,\\
(19) & (C1)1^{k}  & \vdash(k)  & 0^k(C1).
\end{array}$

\smallskip

\noindent Using (5), (19) and (13), we get

$\begin{array}{clll}
(20) & (F0)01^{k+1}0 & \vdash(k + 5) & 10^{k+2}(C0).
\end{array}$

\smallskip

\noindent Using (16), (4), (17) and (9), we get

$\begin{array}{clll}
(21) & 0^2(F0)(10)^k11 & \vdash(4k + 6) & (F0)(10)^{k+2}.
\end{array}$

\smallskip

\noindent Using (16), (20), (15), (17), (9) and (21), we get

$\begin{array}{clll}
(22) & 0^4(F0)(10)^k010 & \vdash(8k + 20) & (F0)(10)^{k+3}1.
\end{array}$

\smallskip

\noindent Using (22) and (21), we get

$\begin{array}{clll}
(23) & 0^6(F0)(10)^k(01)^2 & \vdash(12k + 38) & (F0)(10)^{k+5}.\\
\end{array}$

\smallskip

\noindent By induction on $n$, using (23), we get

$\begin{array}{clll}
(24) & 0^{6n}(F0)(10)^k(01)^{2n} & \vdash(30n^2 + 12kn + 8n) & (F0)(10)^{5n+k}.
\end{array}$

\smallskip

\noindent Using (24), with $k = 2$, we get

$\begin{array}{clll}
(25) & 0^{6k}(F0)(10)^2(01)^{2k} & \vdash(30k^2 + 32k) & (F0)(10)^{5k+2}.
\end{array}$

\smallskip

\noindent Using (16), (2), (4), (11), (18), (8) and (21), we get

$\begin{array}{clll}
(26) & 0^5(F0)(10)^{k+1}0^2 & \vdash(4k + 25) & (F0)(10)^2(01)^k00101.
\end{array}$

\smallskip

\noindent Using (26) and (25), we get

$\begin{array}{clll}
(27) & 0^{6k+5}(F0)(10)^{2k+1}0^2 & \vdash(30k^2 + 40k + 25) & (F0)(10)^{5k+2}00101.
\end{array}$

\smallskip

\noindent Using (25) and (22), we get

$\begin{array}{clll}
(28) & 0^{6k+4}(F0)(10)^2(01)^{2k+1}0 & \vdash(30k^2 + 72k + 36) & (F0)(10)^{5k+5}1.
\end{array}$

\smallskip

\noindent Using (26) and (28), we get

$\begin{array}{clll}
 & 0^{6k+9}(F0)(10)^{2k+2}0^2 & \vdash(30k^2 + 80k + 65) & (F0)(10)^{5k+7}1,
\end{array}$

\noindent and the result is still true for $k = -1$, so we have

$\begin{array}{clll}
(29) & 0^{6k+3}(F0)(10)^{2k}0^2 & \vdash(30k^2 + 20k + 15) & (F0)(10)^{5k+2}1.
\end{array}$

\smallskip

\noindent Using (16), (20), (14), (4), (11), (18), (8) and (21), we get

$\begin{array}{clll}
(30) & 0^5(F0)(10)^k01^20 & \vdash(4k + 23) & (F0)(10)^2(01)^k00101.
\end{array}$

\smallskip

\noindent Using (30) and (25), we get

$\begin{array}{clll}
(31) & 0^{6k+5}(F0)(10)^{2k}01^20 & \vdash(30k^2 + 40k + 23) & (F0)(10)^{5k+2}00101.
\end{array}$

\smallskip

\noindent Using (30) and (28), we get

$\begin{array}{clll}
(32) & 0^{6k+9}(F0)(10)^{2k+1}01^20 & \vdash(30k^2 + 80k + 63) & (F0)(10)^{5k+7}1.
\end{array}$

\smallskip

\noindent Using (16), (20), (14), (4), (17) and (9), we get

$\begin{array}{clll}
(33) & 0^2(F0)(10)^k01^30 & \vdash(4k + 18) & (F0)(10)^{k+3}1.
\end{array}$

\smallskip

\noindent Using (16), (20), (14), (4), (9) and (16), we get

$\begin{array}{clll}
(34) & (F0)(10)^k01^40 & \vdash(2k + 21) & (10)^k1(10)^2(F0)1.
\end{array}$

\smallskip

\noindent Using (34), (4), (17), (11), (18), (8) and (21), we get

$\begin{array}{clll}
(35) & 0^5(F0)(10)^k01^401 & \vdash(4k + 39) & (F0)(10)^2(01)^k00(10)^3.
\end{array}$

\smallskip

\noindent Using (35) and (25), we get

$\begin{array}{clll}
(36) & 0^{6k+5}(F0)(10)^{2k}01^401 & \vdash(30k^2 + 40k + 39) & (F0)(10)^{5k+2}0^2(10)^3.
\end{array}$

\smallskip

\noindent Using (35), (25) and (22), we get

$\begin{array}{clll}
(37) & 0^{6k+9}(F0)(10)^{2k+1}01^401 & \vdash(30k^2 + 80k + 79) & (F0)(10)^{5k+9}.
\end{array}$

\smallskip

\noindent Using (16), (20), (14), (4), (9) and (16), we get

$\begin{array}{clll}
(38) & (F0)(01)^k01^50 & \vdash(2k + 22) & (10)^{k+3}(F0)1.
\end{array}$

\smallskip

\noindent Using (38), (4), (17) and (9), we get

$\begin{array}{clll}
(39) & 0^2(F0)(10)^k01^501 & \vdash(4k + 34) & (F0)(10)^{k+5}.
\end{array}$

\smallskip

\noindent Using (34), (3), (2), (4), (11), (18) and (7), we get

$\begin{array}{clll}
(40) & (F0)(10)^k01^40^4 & \vdash(2k + 43) & (10)^k1(H0)(10)^2(01)^2.
\end{array}$

\smallskip

\noindent Results (1), (21), (29), (27), (22), (31), (32), (33), (36), (37), (39) and (40)
give results (a)--(l) of the theorem.\qed

\smallskip

\noindent Note that the rules (a)--(l) are written in their order of occurrence in the
computation of Turing machine $M_6$ on the blank tape.

\smallskip

Using the rules of this theorem, we have, in 3346 transitions,
$$^\omega0(A0)0^\omega \vdash(6)\quad C(0,15) \vdash(6)\quad C(2,3) \vdash(\ ) \cdots
\vdash(\ ) \quad \mbox{END}$$
We have
$$^\omega0(A0)0\mbox{R(bin($p$))}0^\omega \vdash(6)\quad 
^\omega0(F0)1^40\mbox{R(bin($p$))}0^\omega$$ 
$$= C(0,32p + 15) \vdash (6)\quad C(2,8p + 3) \vdash(14)\quad C(4,2p),$$
so the behavior of Turing machine $M_6$ on inputs 00$x$, $x \in \{0,1\}^*$,
depends on the behavior of configurations $C(n,p)$, and the halting problem for Turing
machine $M_6$ depends on this behavior.

Let us give some informal comments on how Turing machine $M_6$ works.
First, as in the case of the previous machines, machine $M_6$ nibbles
strings $(01)^{2n}$ and transforms them into strings $(01)^{5n}$
(see transitions (23) and (24) in the proof of Theorem \ref{Thm8.1}).
This explains why some transitions of Theorem \ref{Thm8.1} looks like a Collatz-like
behavior of type $2\to 5$.
Second, unlike the previous machines, edge effects are more chaotic.
A study of Turing machine $M_6$ beyond the rules stated in Theorem \ref{Thm8.1}
seems to lead to an infinite set of rules, but we could not prove this
assertion.

\section{Configurations provably stopping}\label{Sect9}

Let $M_7$ be the $6 \times 2$ Turing machine defined by Table~\ref{m7}.

\begin{table}[t]\centering
\begin{tabular}{|c|c|c|}
\hline
$M_7$ &  0  &  1\\
\hline
$A$   & 1RB & 0LB\\
\hline
$B$   & 0RC & 1LB\\
\hline
$C$   & 1RD & 0LA\\
\hline
$D$   & 1LE & 1LF\\
\hline
$E$   & 1LA & 0LD\\
\hline
$F$   & 1RH & 1LE\\
\hline
\end{tabular}
\caption{\label{m7}Machine $M_7$ discovered in October 2000 by Marxen and Buntrock}
\end{table}

We have $s(M_7) > 6.1 \times 10^{925}$ and $\sigma(M_7) > 6.4 \times 10^{462}$.

This machine was discovered in October 2000 by Heiner Marxen and J\"urgen Buntrock.
It was the champion for the busy beaver competition for $6 \times 2$
machines from October 2000 to March 2001.

The following theorem was initially obtained by Munafo \cite{Mu2}.
It gives the rules observed by Turing machine $M_7$.

\begin{thm}\label{Thm9.1}
We have the following transitions between configurations of Turing machine $M_7$.
Let $C(n)$ = $^\omega01^n(B0)0^\omega$. Then

$\begin{array}{clll}
(a) &  ^\omega0(A0)0^\omega & \vdash(1) & C(1),
\end{array}$

and, for all $k \ge 0$,

$\begin{array}{clll}
(b) & C(3k)      & \vdash(54 \times 4^{k+1} - 27 \times 2^{k+3} + 26k + 86)       & C(9 \times 2^{k+1} - 8),\\
(c) & C(3k + 1)  & \vdash(2048 \times(4^k - 1)/3 - 3 \times 2^{k+7} + 26k + 792) & C(2^{k+5} - 8),\\
(d) & C(3k + 2)  & \vdash(3k+8)                                                 &  ^\omega01(H1)(011)^k01010^\omega.
\end{array}$
\end{thm}

\proof
A direct inspection of the transition table gives

$\begin{array}{clll}
(1)  & 11(B0)00   & \vdash(6) & (D1)0101,\\
(2)  & 00(B0)01   & \vdash(8) & (B0)1^4,\\
(3)  & 11(B0)01   & \vdash(6) & (B1)01^3,\\
(4)  & (B0)1      & \vdash(3) & 1(B0),\\
(5)  & 0(B1)      & \vdash(1) & (B0)1,\\
(6)  & 1(B1)      & \vdash(1) & (B1)1,\\
(7)  & 0(D1)      & \vdash(2) & 1(H1),\\
(8)  & 0^31(D1)   & \vdash(6) & (B0)1^4,\\
(9)  & 0^41^2(D1) & \vdash(8) & (B0)1^301^2,\\
(10) & 1^3(D1)    & \vdash(3) & (D1)011.
\end{array}$

\smallskip

\noindent Iterating, respectively, (4), (6) and (10) gives

$\begin{array}{clll}
(11) & (B0)1^k   & \vdash(3k) & 1^k(B0),\\
(12) & 1^k(B1)   & \vdash(k)  & (B1)1^k,\\
(13) & 1^{3k}(D1) & \vdash(3k) & (D1)(011)^k.
\end{array}$

\smallskip

\noindent Using (3), (12), (5) and (11), we get

$\begin{array}{clll}
 & 01^{k+2}(B0)01 & \vdash(4k + 10) & 1^{k+1}(B0)01^3,
\end{array}$

\noindent and the result is still true for $k = -1$, so we have

$\begin{array}{clll}
(14) & 01^{k+1}(B0)01 & \vdash(4k + 6) & 1^k(B0)01^3.
\end{array}$

\smallskip

\noindent By induction on $k$, using (14), we get

$\begin{array}{clll}
(15) & 0^k1^k(B0)01 & \vdash(2k^2 + 4k) & (B0)01^{2k+1}.
\end{array}$

\smallskip

\noindent Using (1), (13) and (7), we get

$\begin{array}{clll}
(16) & 01^{3k+2}(B0)0^2 & \vdash(3k + 8) & 1(H1)(011)^k0101.
\end{array}$

\smallskip

\noindent Using (1), (13) and (9), we get

$\begin{array}{clll}
 & 0^41^{3k+4}(B0)0^2 & \vdash(3k + 14) & (B0)1^3(011)^{k+1}0101,
\end{array}$

\noindent and the result is still true for $k = -1$, so we have

$\begin{array}{clll}
(17) & 0^41^{3k+1}(B0)0^2 & \vdash(3k + 11) & (B0)1^3(011)^k0101.
\end{array}$

\smallskip

\noindent Using (1), (13) and (8), we get

$\begin{array}{clll}
(18) & 0^31^{3k+3}(B0)0^2 & \vdash(3k + 12) & (B0)1^4(011)^k0101.
\end{array}$

\smallskip

\noindent Using (11), (15) and (2), we get

$\begin{array}{clll}
(19) & 0^{k+2}(B0)1^k01 & \vdash(2k^2 + 7k + 8) & (B0)1^{2k+4}.
\end{array}$

\smallskip

\noindent By induction on $k$, using (19), we get

$\begin{array}{clll}
(20) & 0^{2^k(n+5)-5-n-3k}(B0)1^n(011)^k & \vdash(T) & (B0)1^{2^k(n+5)-5},
\end{array}$

\noindent with $T = 2(n+5)^2(4^k-1)/3 - 13(n+5)(2^k-1) + 23k$.

\smallskip

\noindent Using (20), for $n=3$ and $n=4$, we get respectively

$\begin{array}{clll}
(21) & 0^{2^{k+3}-3k-8}(B0)1^3(011)^k & \vdash(128(4^k-1)/3 - 13\times 2^{k+3} + 23k + 104) & (B0)1^{2^{k+3}-5},\\
(22) & 0^{9\times 2^k-3k-9}(B0)1^4(011)^k & \vdash(54\times 4^k - 117\times 2^k + 23k + 63) & (B0)1^{9\times2^k-5}.
\end{array}$

\smallskip

\noindent Using (11), (15), (2), (11), (15), (2) and (11), we get

$\begin{array}{clll}
(23) & 0^{3k+8}(B0)1^k0101 & \vdash(10k^2 + 65k + 112) & 1^{4k+12}(B0).
\end{array}$

\smallskip

\noindent Using (17), (21) and (23), we get

$\begin{array}{clll}
(24) & 0^{2^{k+5}-3k-11}1^{3k+1}(B0)0^2 & \vdash(2048\times(4^k-1)/3 - 3\times 2^{k+7} + 26k + 792) & 1^{2^{k+5}-8}(B0).
\end{array}$

\smallskip

\noindent Using (18), (22) and (23), we get

$\begin{array}{clll}
 & 0^{9\times 2^{k+2}-3k-13}1^{3k+3}(B0)0^2 & \vdash(54\times 4^{k+2} - 27\times 2^{k+4} + 26k + 112) & 1^{9\times 2^{k+2}-8}(B0),
\end{array}$

\noindent and the result is still true for $k = -1$, so we have

$\begin{array}{clll}
(25) & 0^{9\times 2^{k+1}-3k-10}1^{3k}(B0)0^2 & \vdash(54\times 4^{k+1} - 27\times 2^{k+3} + 26k + 86) & 1^{9\times 2^{k+1}-8}(B0).
\end{array}$

\smallskip

\noindent Results (25), (24) and (16) give results (b)--(d) of the theorem.\qed

\smallskip

\noindent Using the rules of this theorem, we have

$$^\omega0(A0)0^\omega \vdash(1)\quad C(1) \vdash(408)\quad C(24) \vdash(14100774)$$
$$C(4600) \vdash(T)\quad C(2^{1538}-8) \vdash(2^{1538}-2)\quad ^\omega01(H1)(011)^p01010^\omega,$$
with $T = 2048 \times (4^{1533}-1)/3 - 3 \times 2^{1540} + 40650$ and $p = (2^{1538}-10)/3$.

So the total time is $s(M_7) = 2048 \times (4^{1533}-1)/3 - 11 \times 2^{1538} + 14141831$,
and the final number of symbols 1 is $\sigma(M_7) = 2 \times (2^{1538} - 10)/3 + 4$.

\smallskip

As in the case of Turing machine $M_4$, we find two ideas that explain how
Turing machine $M_7$ works. First, machine $M_7$ nibbles strings of symbols 0
one by one, transforming them into strings two times longer (see transitions
(14) and (15) in the proof of Theorem \ref{Thm9.1}).
Second, machine $M_7$ nibbles strings $(011)^k$ and, each time a string $011$ is nibbled,
the length of a string of symbols 1 is doubled (see transitions (19) and (20)
in the proof of Theorem \ref{Thm9.1}). This explains why powers of two are obtained.

\smallskip

Note that
$$C(6k+1) \vdash(\ )\quad C(3m) \vdash(\ )\quad C(6p+4) \vdash(\ )\quad C(3q+2) \vdash(\ )\quad
\mbox{END},$$
with $m = (2^{2k+5} - 8)/3$, $p = 3 \times 2^m - 2$, $q = (2^{2p+6} - 10)/3$.
So all configurations $C(n)$ lead to a halting configuration.
Those taking the most time are $C(6k + 1)$.
For example, $C(7) \vdash(t)$ END with $t > 10^{3.9 \times 10^{12}}$.
More generally, $C(6k + 1) \vdash(t(k))$ END with $t(k) > 10^{10^{10^{(3k+2)/5}}}$.

\section{Conclusion}

We discuss two questions as a conclusion to this article.

\smallskip

\noindent {\bf A.} \emph{How simulating Collatz-like functions allows Turing machines to
achieve high scores?}

Lagarias \cite{La85} noted that the successive iterates of the 3$x$+1 function $T$
have an irregular behavior. For example, 7 iterations of function $T$ on $n = 26$
lead to the value 1, but 70 iterations are necessary on $n = 27$.
It seems that many Collatz-like functions have the same irregular behavior.
Iterating them on small numbers may produce very long runs before stopping.

Adding parameters may increase the number of iterations by allowing the iterated
values to range the set of parameters before stopping. The pure Collatz-like
function with parameter $g_3(n,p)$ presented in Section \ref{Sect5} is particularly
illustrative.

Another way to high scores is given by exponential Collatz-like functions such
as function $g_4$ in Section \ref{Sect6}. Only five iterations are performed on a blank tape,
but exponential growth ensures a high score.

Irregular behavior is a condition for a Collatz-like function to be eligible to the
busy beaver competition. Another condition is, of course, being computable by a very
small Turing machine.

\smallskip

\noindent {\bf B.} \emph{Are some universal devices more natural than others?}

Conway \cite{Co72} proved that there is no algorithm that, given as inputs
a Collatz-like function $g$ and two integers $n$, $p$, outputs an answer yes or no
to the question: Does there exist a positive integer $k$ such that $g^k(n) = p$?
Conway \cite{Co72, Co87} also proved that Collatz-like functions can be used to
simulate all computable (also called recursive) functions.
These properties can be summed up by writing that Collatz-like functions
provide a universal model of computation with a m-complete decision problem.

Many universal models of computation are known: Turing machines, tag-systems,
cellular automata, Diophantine equations, etc. (see \cite{MM10}).
Of course, any universal model can simulate and be simulated by any other
universal model. But it is Collatz-like functions, and not another model,
that appear naturally in this study. Their unexpectedly pervasive presence
leads to wonder about the significance of their status among mathematical beings.

\section*{Acknowlegement}
We thank an anonymous referee for many helpful suggestions.

\end{document}